\documentclass[review]{elsarticle}

\usepackage{hyperref}
\usepackage{graphicx}

\usepackage{array,color}
\usepackage{amsbsy}
\usepackage{amsmath}
\usepackage{amssymb}
\usepackage{amsfonts}
\usepackage{epsf}
\usepackage{epsfig}
\usepackage{graphicx}
\usepackage{fancybox,fancyhdr,times,subfigure}
\usepackage{epsfig, color, rotate, colordvi}
\usepackage{a4wide}
\usepackage{natbib}
\usepackage{mathrsfs}
\usepackage{bm}
\usepackage{multirow}
\usepackage{pgfplots}
\usepackage{diagbox}
\usepackage[titletoc]{appendix}
\usepackage{titlesec}
\usepackage{algorithm}
\usepackage{algorithmic}

\newcommand{\bbx}{{\bf x}}

\newcommand{\bbX}{{\bf X}}

\newcommand{\bbI}{{\bf I}}

\newcommand{\bqn}{\begin{eqnarray*}}
\newcommand{\eqn}{\end{eqnarray*}}
\newcommand{\bqa}{\begin{eqnarray}}
\newcommand{\eqa}{\end{eqnarray}}

\newcommand{\bqnn}{\begin{eqnarray}}
\newcommand{\eqnn}{\end{eqnarray}}

\newcommand{\non}{\nonumber\\}

\newtheorem{thm}{Theorem}
\newtheorem{corollary}{Corollary}
\newtheorem{proof}{Proof}
\newtheorem{lem}[thm]{Lemma}
\newtheorem{rmk}{Remark}
\newtheorem{condition}{Condition}

\begin{document}

\begin{frontmatter}

\title{Application of fused graphical lasso to statistical inference for multiple sparse precision matrices}

\author[mymainaddress]{Qiuyan Zhang}
%\ead{qyzhang@cueb.edu.cn}
\address[mymainaddress]{School of Statistics, Capital University of Economics and Business, Beijing}

\author[mysecondaryaddress]{Zhidong Bai}
\address[mysecondaryaddress]{KLASMOE and School of Mathematics and Statistics, Northeast Normal University, Changchun, Jilin}

\author[mymainaddress]{Lingrui Li}

\author[add3]{Hu Yang\corref{mycorrespondingauthor}}
\cortext[mycorrespondingauthor]{Corresponding author}
\ead{hu.yang@cufe.edu.cn}
\address[add3]{School of Information, Central University of Finance and Economics, Beijing}

\begin{abstract}
 In this paper, the fused graphical lasso (FGL) method is used to estimate multiple precision matrices from multiple populations simultaneously.
 The lasso penalty in the FGL model is a restraint on sparsity of precision matrices, and a moderate penalty on the two precision matrices from distinct groups restrains the similar structure across multiple groups.
 In high-dimensional settings, an oracle inequality is provided for FGL estimators, which is necessary to establish the central limit law.
 We not only focus on point estimation of a precision matrix, but also work on hypothesis testing for a linear combination of the entries of multiple precision matrices.
 Inspired by Jankov$\acute{a}$ and van de Geer [confidence intervals for high-dimensional inverse covariance estimation, {\it Electron. J. Stat.} {\bf 9}(1) (2015) 1205-1229.], who investigated a de-biasing technology to obtain a new consistent estimator with known distribution for implementing the statistical inference, we extend the statistical inference problem to multiple populations, and propose the de-biasing FGL estimators.
 The corresponding asymptotic property of de-biasing FGL estimators is provided.
 A simulation study shows that the proposed test works well in high-dimensional situations.
\end{abstract}


\end{frontmatter}

\section{Introduction}
 Undirected graphical models are popular tools for representing the network structure of data and have been widely applied in many domains, such as machine learning, genetics, and finance.
 Letting $\bbx=(\bbx^{1},...,\bbx^{p})^{T}$ be a p-variate normal random vector with mean vector $\mu$ and covariance $\Sigma_0$ ($\Sigma_0$ is positive definite),
 the precision matrix (or concentration matrix) is denoted the inverse of the covariance matrix, i.e., $\Theta_0:=\Sigma_0^{-1}$.
 The graphical models capture conditional dependence relationships between random variables via non-zero entries in a precision matrix.
 If $\Theta_{0ij}\neq 0$, $\bbx^{i}$ and $\bbx^{j}, i,j=1,...,p$ are dependent on each other, given all other variables.
 Meanwhile, the zero entries in the precision matrix correspond to pairs of variables that are conditionally independent given other variables.
 Therefore, the graph model is closely related to the precision matrix.
 Estimating and testing of a precision matrix have been a rapidly growing research direction in the past few years.

 Letting $\bbx_1,...,\bbx_n$ be a sequence of independent and identically distributed (i.i.d.) observations from the population $\bbx$, $\bbX_{p\times n}:=(\bbx_1,...,\bbx_n)$.
 A natural estimator of the precision matrix is the inverse of the sample covariance matrix $\widehat\Sigma$, where $\widehat\Sigma=\frac{1}{n}\bbX^{T}\bbX$.
 On one hand, in high-dimensional settings, Johnstone \cite{johnstone2001distribution} proposed that the eigenvalues of the sample covariance matrix do not converge to the corresponding eigenvalue of the population covariance matrix for $\Sigma=\bbI$.
 Consequently, this estimator becomes invalid when the dimension $p$ is comparable to the sample size $n$.
 On the other hand, the sample covariance matrix is singular in a $p>n-1$ setting.
 This will produce non-negligible errors in using $\widehat\Sigma_n^{-1}$ to estimate $\Theta_0$.
 In addition, a sparse (i.e., many entries are either zero or nearly so) assumption for a high-dimensional precision matrix is essential, since the zero entries imply the conditional independence structures, which are what we are most concerned with in the graphical model.
 In general, $\widehat\Sigma^{-1}$ does not have a sparsity construction.
 How to estimate the sparse precision matrix in high-dimensional settings is an intractable problem.

 In recent years, various proposals have been put forward for estimating a precision matrix in high-dimensional situations, among which the graphical model with sparsity-promoting penalties is valid for obtaining a sparse estimator.
 By applying the $l_1$ (lasso penalty) to the entries of the concentration matrix, Yuan and Lin \cite{yuan2007model} proposed a max-det algorithm to obtain the estimator of $\Theta_0$. The convergence result of the estimator is derived under a $p$ fixed assumption.
 Using a coordinate descent procedure, Friedman et al. \cite{friedman2008sparse} provided an algorithm for solving a graphical Lasso estimator that is remarkably fast, even if $p>n$.
 Rothman et al. \cite{Rothman2008Sparse} investigated a sparse permutation invariant covariance estimator, and
 established a convergence rate of the estimator in the Frobenius norm as both data dimension $p$ and sample size $n$ are allowed to grow, and showed that the rate explicitly depends on how sparse the true concentration matrix is.
 For additional theoretical details on penalized likelihood methods for graphical models, see Fan et al. \cite{fan2009network}, Ravikumar et al. \cite{ravikumar2011high}, Xue and Zou \cite{xue2012regularized}, and Yuan et al.\cite{yuan2019constrained}.

 The above-mentioned methods focus on estimating a single graphical model, but joint estimators better recover the truth graphs compared with separate estimations for data from multiple graphical models sharing the similarities structure with other, but not identical models.
 Guo et al. \cite{guo2011joint} studied joint estimation of precision matrices that have a hierarchical structure assumption.
 Liu and Lee proposed a joint estimator of multiple precision matrices under an assumption that precision matrices decompose into the sum of two components.
 A fused graphical lasso was proposed by Danaher et al. \cite{Danaher2014The} with a penalty imposing a similar structure of a precision matrix across groups.
 Supposing that $\bbX^{[k]}_{p\times n_k}:=(\bbx_1^{[k]},...,\bbx_{n_k}^{[k]})$ are sample matrices, and $\bbx_{i}^{[k]}\in R^{p} (i=1,...,n_k)$ are sampled i.i.d. from a distribution with mean $\mu^{[k]}$ and covariance $\Sigma_0^{[k]}$, for $k=1,...,K$, we assume $\mu^{[k]}=0$ without loss of generality.
 To simplify notation, we omit the subscript of $\bbX^{[k]}_{p\times n_k}$, and denote the sample matrices as $\bbX^{[k]}$.
 The population precision matrix is defined as the inverse of the population covariance matrix, i.e., $\Theta_0^{[k]}=(\Sigma_0^{[k]})^{-1}$.
 The estimators of precision matrices $\{\Theta_0^{[k]}\}$ are investigated by minimizing the negative penalized log likelihood
\bqa\label{op}
\{\widehat{\Theta}^{[k]}\}={\arg\min}_{\{\Theta^{[k]}\in\mathcal{S}^{++}\}} \sum_{k} \{tr(\widehat\Sigma^{[k]}\Theta^{[k]})-\log \det(\Theta^{[k]})\}+{\bf{P}}(\{\Theta^{[k]}\}),
\eqa
 where ${\bf{P}}(\{\Theta^{[k]}\})$ denotes the penalty function, the $\{\widehat{\Theta}^{[k]}\}$ are the minimizers of (\ref{op}), and we optimize over the symmetric positive-definite matrices set $\mathcal{S}^{++}$.
 The fused graphical lasso (FGL) is the solution to optimization problem (\ref{op}) with the fused lasso penalty
\bqa\label{Genpen}
{\bf{P}}(\{\Theta^{[k]}\})=\lambda\sum_{k=1}^{K}||(\Theta^{[k]})^{-}||_{1}
+\rho\sum_{k<k'}||(\Theta^{[k]}
-\Theta^{[k']})^{-}||_1,
\eqa
 where $\lambda$ and $\rho$ are non-negative regularization parameters, $(\Theta^{[k]})^{-}$ represents the matrix obtained by setting the diagonal elements of $(\Theta^{[k]})$ to zero, and $||\cdot||_1$ denotes the $l_1$ norm of a vector or matrix.
 It is reasonable to restrict non-diagonal elements of $\Theta^{[k]}$, since we are most concerned with the conditional independence cross-different variables.
 Note that the first term in (\ref{Genpen}) is the classical lasso penalty, which shrinks the coefficients toward $0$ as $\lambda$ increases.
 It guarantees discovery of the sparse estimators $\{\widehat{\Theta}^{[k]}\}$ of the model.
 The penalty on $(\Theta^{[k]}
-\Theta^{[k']})^{-}$ indicates that the elements of $\widehat{\Theta}^{[1]},...,\widehat{\Theta}^{[K]}$ have a similar network structure across classes.

 An approach for the estimation of the joint graphical models largely relies on penalized estimation.
 The penalty biases the estimates toward the assumed structure, which makes hypothesis tests for precision matrices more challenging.
 Work on statistical inference for low-dimensional parameters in graphical models has recently been carried out (%Barber and Kolar \cite{foygel2015rocket};
Jankov$\acute{a}$ and van de Geer \cite{jankova2015confidence}; Jankov$\acute{a}$ and van de Geer \cite{jankova2017honest}; Ren et al. \cite{ren2015asymptotic}; Yu et al. \cite{yu2020simultaneous}) based on the $l_1$-penalized estimator.
 Jankov$\acute{a}$ and van de Geer \cite{jankova2015confidence} provided a de-biasing technique to obtain a new consistent estimator with known distribution.
 However, these approaches were developed only in the setting in which the parameters of one graph are inferred.
 In contrast, studies of inference techniques using estimators obtained from cross-group penalization are much fewer.
 The work on statistical inference for multiple graphical models is an interesting area open for future research.
 Inspired by Jankov$\acute{a}$ and van de Geer \cite{jankova2015confidence}, we not only give FGL estimators of multiple precision matrices from co-movement data, but also test the linear combination of the entries of these precision matrices.
 The core of the proposed method is based on the de-biasing technique, and we implement statistical inference of the precision matrices under high-dimensional settings according to the proposed central limit theorem.

 The rest of this paper is organized as follows.
 In Section 2, we give the oracle inequality for multiple estimators with a FGL penalty and its weighted version.
 Testing the hypothesis for the linear combination of corresponding entries of multiple precision matrices is also considered in this section.
 Based on de-biasing technology, the CLT of the proposed statistics for multiple populations is also derived in Section 2.
 In Section 3, we report the results of  simulations.
 All technical details are relegated to the Appendix.

\section{Main results}

 We assume following notation throughout the paper.
 For a matrix $A=(a_{ij})_{i,j=1}^{p}$, we denote $(A)_{ij}$ its $(i,j)$-entry, or denote its $(i,j)$-entry as $A_{ij}$ to simplify the notation.
 We write $|A|$ for the determinant of $A$, and the trace of matrix $A$ is denoted $tr(A)$.
 Letting $A^{+}=diag(A)$ for a diagonal matrix with the same diagonal as $A$, $A^{-}=A-A^{+}$.
 $||A||_F^2=\sum_{i,j}a_{ij}^2$ denotes the Frobenius norm (also known as the matrix 2-norm).
 We use the notation $||A||_{\infty}=\max_{i,j}|a_{ij}|$ for the supremum norm of a matrix $A$, and $|||A|||_1:=\max_{j}\sum_{i}|a_{ij}|$ for the $l_1$-operator norm.

 We write $f(n)=\mathcal{O}(g(n))$ if $f(n)\leq cg(n)$ for some constant $c<\infty$, and $f(n)=\Omega(g(n))$ if $f(n)\geq c'g(n)$ for some constant $c'>0$.
 The notation $f(n)\asymp g(n)$ means that $f(n) = \mathcal{O}(g(n))$ and $f(n)=\Omega(g(n))$.
 In the common high-dimensional setting, the dimension $p$ is allowed to grow to infinity.
 The dimension is comparable, substantially larger or smaller than the sample size.
 We set sample sizes $ n_1\asymp ...\asymp n_K \asymp n$ throughout the paper, and $n^*=n_1+...+n_K$ going to infinity.
 Furthermore, for notational simplicity, we assume that $n_1=...=n_K=n$.

\subsection{Oracle inequality}

 To obtain the oracle inequality of multiple estimators of FGL models, we introduce some notation related to the sparsity assumptions on the entries of the true precision matrix.
 Letting
\bqa
S_k:=\{(i,j):\Theta^{[k]}_{0ij}\neq 0,i\neq j\},
\eqa
 where $\Theta^{[k]}_{0ij}$ is the $(i,j)$-entry of $\Theta^{[k]}_{0}$ and $s_k=|S_k|$ is the cardinality of $S_k$, we adopt the boundedness of the eigenvalues of the true precision matrix and certain tail conditions proposed by Jankov$\acute{a}$ and Van De Geer \cite{jankova2015confidence}.
\begin{condition}[Bounded eigenvalues]
There exist universal constants $L$ for $k$ such that
\begin{equation*}
0 <L < \Lambda_{\mathrm{min}}(\Theta_{0}^{[k]}) \leq \Lambda_{\mathrm{max}}(\Theta_{0}^{[k]}) <1/L<\infty,
\end{equation*}
\end{condition}
 where $\Lambda_{\mathrm{min}}$ and $\Lambda_{\mathrm{max}}$ denote the minimum and maximum eigenvalues of a matrix, respectively.

\begin{condition}[Sub-Gaussianity vector condition]
The observations $\bbx_i^{[k]}$, $i=1,\dots,n_k$, are uniformly sub-Gaussian vectors in the respective groups.
\end{condition}

 We propose the oracle inequality for FGL lasso under the $K=2$ situation.

\begin{thm}\label{Orcalebounds}
 Supposing that Conditions 1 and 2 hold, for $k=1,2$,
 tuning parameter $\lambda$ satisfying $2(\rho+\lambda_0)\leq\lambda\leq c/8L$, and $\frac{8\lambda^2(s_1+s_2)}{c}+\frac{4p\lambda_0^2}{c}\leq \lambda_0 /2L$. On the set $\{\max_{k}||\widehat{\Sigma}^{[k]}-\Sigma_{0}^{[k]}||_{\infty}\leq\lambda_0\}$, $k=1,2$,
 it holds that
\begin{equation}
\begin{split}
c\sum_{k=1}^{2}||\widehat{\Theta}^{[k]}-\Theta_{0}^{[k]}||_{F}^{2}
+\lambda\sum_{k=1}^{2}||(\widehat{\Theta}^{[k]}-\Theta_{0}^{[k]})^{-}||_1
\leq \frac{8\lambda^2(s_1+s_2)}{c}
+\frac{4p\lambda_0^2}{c},
\end{split}
\end{equation}
 and
\begin{equation}
\begin{split}
\sum_{k=1}^2|||\widehat{\Theta}^{[k]}-\Theta_{0}^{[k]}|||_1\leq\frac{4\lambda(8s_1+8s_2+p)}{c},
\end{split}
\end{equation}
where $c=1/(8L^2)$.
\end{thm}

\begin{rmk}
 From the inequality, we must select $\lambda$ so that $\lambda p\to 0$ as $n\to \infty$ to ensure consistency, which is not satisfied by a sub-Gaussianity random vector.
 Thus, the condition $\lambda p\to 0$ excludes the $p\gg n$ situation.
\end{rmk}

 The FGL does not take into account that the variables have, in general, different scaling.
 Thus, we consider the weighted FGL.
 The minimizer of the optimization problem (\ref{op}) with weighted FGL penalty
\bqa\label{wei}
{\bf{P}}(\{\Theta^{[k]}\})=\lambda\sum_{k}
\sum_{i\neq j}\widehat{W}_{ii}^{[k]}\widehat{W}_{jj}^{[k]}|\Theta_{ij}^{[k]}|
+\rho\sum_{k<k'}\sum_{i\neq j}|\widehat{W}_{ii}^{[k]}\widehat{W}_{jj}^{[k]}\Theta_{ij}^{[k]}
-\widehat{W}_{ii}^{[k']}\widehat{W}_{jj}^{[k']}\Theta_{ij}^{[k']}|
\eqa
 is denoted $\{\widehat{\Theta}_w^{[k]}\}$, where $\widehat{W}^{[k]}=\big[diag(\widehat{\Sigma}^{[k]})\big]^{\frac12}$.
 Further, the population correlation matrix is denoted $R_0^{[k]}$ and the sample correlation matrix is denoted
\bqa
\widehat{R}^{[k]}=(\widehat{W}^{[k]})^{-1}\widehat{\Sigma}^{[k]}
 (\widehat{W}^{[k]})^{-1}.
\eqa
 If we substitute $\widehat{R}^{[k]}$ for $\widehat{\Sigma}^{[k]}$, the minimizer of
\bqa\label{normalized}
{\arg\min}_{\{\Theta^{[k]}\in\mathcal{S}^{++}\}} \sum_{k} \{tr(\widehat{R}^{[k]}\Theta^{[k]})-\log \det(\Theta^{[k]})\}+{\bf{P}}(\{\Theta^{[k]}\})
\eqa
 with a FGL penalty (\ref{Genpen}) is denoted $\{\widehat{\Theta}_R^{[k]}\}$, which is a matter of estimating the parameter by the normalized data.
 Then,
\bqa
\widehat{\Theta}_R^{[k]}=\widehat{W}^{[k]}\widehat{\Theta}_w^{[k]}\widehat{W}^{[k]},
\eqa
 which means, essentially, that $\widehat{\Theta}_R^{[k]}$ are the estimators of $\Theta_{R0}^{[k]}:=\big(R_0^{[k]}\big)^{-1}$.

\begin{thm}\label{weiandr}
 Under the conditions of Theorem \ref{Orcalebounds}, on the set $\{\max_{k}||\widehat{R}^{[k]}-R_{0}^{[k]}||_{\infty}\leq\lambda_0\}$, $k=1,2$, it holds that
\begin{equation}
\begin{split}
c\sum_{k=1}^{2}||\widehat{\Theta}_R^{[k]}-\Theta_{R0}^{[k]}||_{F}^{2}
+\lambda\sum_{k=1}^{2}||(\widehat{\Theta}_R^{[k]}-\Theta_{R0}^{[k]})^{-}||_1
\leq \frac{8\lambda^2(s_1+s_2)}{c},
\end{split}
\end{equation}
 \begin{equation}
\begin{split}
\sum_{k=1}^2|||\widehat{\Theta}_R^{[k]}-\Theta_{R0}^{[k]}|||_1
\leq\frac{32\lambda(s_1+s_2)}{c},
\end{split}
\end{equation}
 and
\begin{equation}
\begin{split}
\sum_{k=1}^2|||\widehat{\Theta}_w^{[k]}-\Theta_{0}^{[k]}|||_1
\leq\frac{32\lambda(s_1+s_2)}{c}.
\end{split}
\end{equation}
\end{thm}

 It is natural to extend this conclusion to the $K>2$ FGL model. For $k=1,...,K$ and the $K>2$ situation, we obtain the following theorem.

\begin{thm}[Multiple FGL model]\label{FGL}
 Supposing that Conditions 1 and 2 hold, for $K>2$, $2\left(\frac{K(K-1)}{2}\rho+\lambda_0\right)\leq\lambda\leq c/8L$, and $\frac{8\lambda^2\sum_{k=1}^{K}s_k}{c}+\frac{2Kp\lambda_0^2}{c}\leq \lambda_0 /2L$, on the set $\{\max_{k}||\widehat{\Sigma}^{[k]}-\Sigma_{0}^{[k]}||_{\infty}\leq\lambda_0\}$, $k=1,...,K$,
 it holds that
\begin{equation}\label{FGLeq1}
\begin{split}
c\sum_{k=1}^{K}||\widehat{\Theta}^{[k]}-\Theta_{0}^{[k]}||_{F}^{2}
+\lambda\sum_{k=1}^{K}||(\widehat{\Theta}^{[k]}-\Theta_{0}^{[k]})^{-}||_1
\leq \frac{8\lambda^2\sum_{k=1}^{K}s_k}{c}+\frac{2Kp\lambda_0^2}{c}
\end{split}
\end{equation}
 and
\begin{equation}\label{FGLeq2}
\begin{split}
\sum_{k=1}^K|||\widehat{\Theta}^{[k]}-\Theta_{0}^{[k]}|||_1\leq
\frac{2K\lambda\left(8\sum_{k=1}^{K}s_k+\frac{Kp}{2}\right)}{c}.
\end{split}
\end{equation}
\end{thm}

\begin{thm}[Multiple FGL model for weighted version]\label{FGLwei}
 Under the conditions of Theorem \ref{FGL}, on the set $\{\max_{k}||\widehat{R}^{[k]}-R_{0}^{[k]}||_{\infty}\leq\lambda_0\}$, $k=1,2$, it holds that
\begin{equation}
\begin{split}
 c\sum_{k=1}^{K}||\widehat{\Theta}_R^{[k]}-\Theta_{R0}^{[k]}||_{F}^{2}
 +\lambda\sum_{k=1}^{K}||(\widehat{\Theta}_R^{[k]}-\Theta_{R0}^{[k]})^{-}||_1
\leq
\frac{8\lambda^2\sum_{k=1}^{K}s_k}{c},
\end{split}
\end{equation}
\begin{equation}
\begin{split}
\sum_{k=1}^{K}|||\widehat{\Theta}_R^{[k]}-\Theta_{R0}^{[k]}|||_{1}
\leq
\frac{16K\lambda\sum_{k=1}^{K}s_k}{c},
\end{split}
\end{equation}
 and
\begin{equation}
\begin{split}
\sum_{k=1}^{K}|||\widehat{\Theta}_w^{[k]}-\Theta_{0}^{[k]}|||_1
\frac{16K\lambda\sum_{k=1}^{K}s_k}{c}.
\end{split}
\end{equation}
\end{thm}

\subsection{Asymptotic property}

 We not only focus on the point estimation of multiple precision matrices, but also on hypothesis testing for the linear combination of the entries of the precision matrices over two groups.
% Let $\theta_{ij}^{[k]}$ be the $(i,j)$-entry of the precision matrix $\Theta_0^{[k]}$.
 One may want to test whether the elements of the precision matrix over two groups are equal:
\bqa\label{hy1}
H_0: \Theta^{[1]}_{0ij}=\Theta^{[2]}_{0ij} \quad vs. \quad H_1: \Theta^{[1]}_{0ij}\neq \Theta^{[2]}_{0ij}.
\eqa

 To test Hypothesis (\ref{hy1}), we aim to obtain confidence intervals for estimators based on the de-biasing technique, which is imposed for eliminating the bias associated with the penalty.
 The de-biasing estimator is defined as $\widehat{{\Theta}}^{[k]}_d=2\widehat{\Theta}^{[k]}-\widehat{\Theta}^{[k]}
\widehat{\Sigma}^{[k]}\widehat{\Theta}^{[k]}$.
 The difference between the de-biasing estimator and the true value can be decomposed into two parts as follows:
\bqa
\widehat{{\Theta}}^{[k]}_d-\Theta_0^{[k]}=\Xi^{[k]}+\Upsilon^{[k]},
\eqa
 where
\bqa
&\Xi^{[k]}=-\Theta_0^{[k]}
(\widehat{\Sigma}^{[k]}-\Sigma^{[k]}_0)\Theta_0^{[k]},\\
&\Upsilon^{[k]}=-(\widehat{\Theta}^{[k]}-\Theta_0^{[k]})
 (\widehat{\Sigma}^{[k]}-\Sigma_0^{[k]})\Theta_0^{[k]}-(\widehat{\Theta}^{[k]}-\Theta_0^{[k]})
 (\widehat{\Sigma}^{[k]}\widehat{\Theta}^{[k]}-\bbI_p).
\eqa
 Under the compatibility conditions, Jankov{\'a} and van~de Geer \cite{jankova2018inference} proposed that the $(i,j)$-entry of $\widehat{{\Theta}}^{[k]}_d-\Theta_0^{[k]}$ has an asymptotic normality property, and $\sqrt{n}||\Upsilon^{[k]}||_{\infty}$ converges to zero in probability.
 Thus, for testing Hypothesis (\ref{hy1}), we construct the testing statistic
\bqa\label{teststa}
T_{ij}:=\left(\widehat{{\Theta}}_{d}^{[1]}-\widehat{{\Theta}}_{d}^{[2]}\right)_{ij}=\left[2\widehat{\Theta}^{[1]}-\widehat{\Theta}^{[1]}
\widehat{\Sigma}^{[1]}\widehat{\Theta}^{[1]}-(2\widehat{\Theta}^{[2]}-\widehat{\Theta}^{[2]}
\widehat{\Sigma}^{[2]}\widehat{\Theta}^{[2]})\right]_{ij}
\eqa
using de-biasing estimators.

For $K=2$, we let
\bqa
s=\max\{s_1,s_2\},\quad d=\max\{d_1,d_2\},
\eqa
 where
\bqa
d_k=\max_{j=1,...,p}|D_{j}^{[k]}|,\quad D_{j}^{[k]}=\{(i,j):\Theta^{[k]}_{0ij}\neq 0,i\neq j\}.
\eqa
 Next, we establish the central limit theorem for $T_{ij}$.

\begin{thm}\label{CLT}
Assuming Conditions 1, 2, and $\lambda\asymp \rho\asymp \sqrt{\log p/n}$ and $(p+s)\sqrt{d}=o(\sqrt{n}/\log p)$,
 it holds that
\bqa\label{linear}
\widehat{{\Theta}}^{[1]}_d-\widehat{{\Theta}}^{[2]}_d-(\Theta_0^{[1]}-\Theta_0^{[2]})=\Xi^{[1]}-\Xi^{[2]}
+rem,
\eqa
 where
\bqa
||rem||_{\infty}=||\Upsilon^{[1]}-\Upsilon^{[2]}||_{\infty}=o_{p}(1/\sqrt{n}),
\eqa
and $o_p$ denotes the convergence in probability. Moreover,
\bqa
\sqrt{n}\big[T_{ij}-{\bf \Theta}_{0ij}\big]\to_{D}N(0,\sigma_{ij}^2),
\eqa
 where ${\bf \Theta}_{0ij}=(\Theta_0^{[1]}-\Theta_0^{[2]})_{ij}$.
\end{thm}

 To complete the testing procedure, we use the consistent estimator
$\hat\sigma_{ij}^2=(\widehat{\Theta}^{[1]})_{ii}(\widehat{\Theta}^{[1]})_{jj}
 +(\widehat{\Theta}^{[1]})_{ij}^2+(\widehat{\Theta}^{[2]})_{ii}(\widehat{\Theta}^{[2]})_{jj}
 +(\widehat{\Theta}^{[2]})_{ij}^2$ for Theorem \ref{CLT}.
 Theorem \ref{CLT} provide a practical and efficient way of obtaining the p value and critical value for the test statistic.
 Under a null hypothesis, we observe that $\Theta_{0ij}^{[1]}-
\Theta_{0ij}^{[2]}=0$.
 For an $\alpha$ level of significance, we reject $H_0$ if $|\sqrt{n}T_{ij}/\hat\sigma_{ij}^2| > \xi_{\alpha/2}$, where $\xi_{\alpha}$ is the $1-\alpha$ upper quantile of the standard normal distribution.

 Theorem \ref{CLT} requires a stronger sparsity condition than the corresponding oracle-type inequality in Theorem \ref{Orcalebounds}.
 According to the convergence rate of $(p+s)\sqrt{d}$, Theorem \ref{CLT} applies to the $p\ll n$ situation.
 For $p\gg n$, we provide the following theorem.

\begin{thm}\label{weiverclt}
 Assuming Conditions 1, 2, and $\lambda\asymp \rho \asymp \sqrt{\log p/n}$ and $s\sqrt{d}=o(\sqrt{n}/\log p)$, for the $p\ll n$ regime, the equation (\ref{linear}) holds with $\widehat{\Theta}_w^{[k]}$, where
\bqa
||rem||_{\infty}=o_{p}(1/\sqrt{n}).
\eqa
 In addition,
\bqa
\sqrt{n}\big[T_{wij}-{\bf \Theta}_{0ij}\big]\to_{D}N(0,\sigma_{ij}^2),
\eqa
 where $T_{wij}=(2\widehat{\Theta}_{w}^{[1]}-\widehat{\Theta}_{w}^{[1]}
\widehat{\Sigma}^{[1]}\widehat{\Theta}_{w}^{[1]})_{ij}-(2\widehat{\Theta}_w^{[2]}
-\widehat{\Theta}_w^{[2]}
\widehat{\Sigma}^{[2]}\widehat{\Theta}_w^{[2]})_{ij}$.
\end{thm}

 We do not need to impose the so-called irrepresentability condition on $\Sigma$ to derive the theoretical properties of our estimators, in contrast to Brownlees et al. \cite{brownlees2018realized}.

 In addition, for the multi-sample precision matrix hypothesis problem, one may want to test
 a linear hypothesis testing problem:
\bqa
H_0: a_1\Theta^{[1]}_{0ij}+...+a_K\Theta^{[K]}_{0ij}=0 \quad vs. \quad H_1: \mbox{not} \ \ H_0,
\eqa
 where $a_1,...,a_K$ are known constants.
 Similar to the two-sample case, we proposed the test statistic
\bqa
a_1\widehat{{\Theta}}^{[1]}_{dij}
+...+a_K\widehat{{\Theta}}^{[K]}_{dij}.
\eqa
 For the $K>2$ multiple situation, we assume $s=\max\{s_1,...,s_K\}$ and $d=\max\{d_1,...,d_K\}$.
 Consequently, we establish the asymptotic normality of the proposed statistic in the following corollary, i.e., Corollary \ref{multipletest}.

\begin{corollary}\label{multipletest}
 Under the assumptions of Theorem \ref{CLT}, it holds that
\bqa\label{linear}
f\big(\widehat{{\Theta}}^{[1]}_{d},...,
\widehat{{\Theta}}^{[K]}_{d}\big)
-f\big(\Theta_{0}^{[1]},...,\Theta_{0}^{[K]}\big)
=f\big(\Xi^{[1]},...,\Xi^{[K]}\big)
+rem,
\eqa
\bqa\label{rem1oversn}
||rem||_{\infty}=||f\big(\Upsilon^{[1]},...,\Upsilon^{[K]}\big)||_{\infty}=o_{p}(1/\sqrt{n}),
\eqa
 where $f(x_1,...,x_K)=a_1x_1+...+a_Kx_K$.
 In addition,
\bqa\label{cltmul}
\sqrt{n}\big[T_{ij}-{\bf \Theta}_{0ij}\big]
\to_{D}N(0,\sigma_{ij}^2),
\eqa
 where $T_{ij}=f\left(\widehat{{\Theta}}^{[1]}_{dij},...,\widehat{{\Theta}}^{[K]}_{dij}\right)$ and ${\bf \Theta}_{0ij}=f\left(\Theta_{0ij}^{[1]},...,\Theta_{0ij}^{[K]}\right)$.
\end{corollary}
 The asymptotic variance $\sigma_{ij}$ in Corollary \ref{multipletest} is unknown, so to construct confidence intervals we use a consistent estimator
\bqa
\hat\sigma_{ij}^2=f_v\big(\big[(\widehat{\Theta}^{[1]})_{ii}(\widehat{\Theta}^{[1]})_{jj}
 +(\widehat{\Theta}^{[1]})_{ij}^2\big],...,\big[(\widehat{\Theta}^{[K]})_{ii}(\widehat{\Theta}^{[K]})_{jj}
 +(\widehat{\Theta}^{[K]})_{ij}^2\big]\big),
\eqa
 where $f_v(x_1,...,x_K)=a_1^2x_1+...+a_K^2x_K$.
 In addition, a weighted version is proposed as follows.

\begin{corollary}
 Under the assumptions of Theorem \ref{weiverclt}, the residual term in (\ref{rem1oversn}) converges in probability with rate $1/\sqrt{n}$, and CLT in (\ref{cltmul}) holds by replacing $\widehat{\Theta}^{[k]}$ by $\widehat{\Theta}_w^{[k]}$, which is obtained by solving the weighted FGL optimization problem.
\end{corollary}

\section{Numerical study}

Simulation experiments were carried out to evaluate the performance of the proposed de-biasing FGL test.
 We considered the sparse graphical model, and a random sample was generated from the multivariate normal distribution $N(0_p,(\Theta_0^{[k]})^{-1})$ with a population covariance matrix defined as the inverse of the population precision matrix.

 To solve the graphical lasso problem with a certain penalty, we refer to the alternating direction method of multiplier (ADMM) algorithm, since it is guaranteed to converge to the global optimum. For more details, the reader is referred to Boyd et al. \cite{boyd2004convex} and Danaher et al. \cite{Danaher2014The}.
 When an objective method for selecting tuning parameters $\lambda$ and $\rho$ is required, the approximations of the Akaike information criterion (AIC), Bayesian information criterion, or cross-validation method can be used to select tuning parameters.
 The AIC method was chosen for the following simulation, and $\lambda$ and $\rho$ both range from $0.05$ to $0.3$ with a step of $0.0086$, where the step is derived by $(0.3-0.05)/(30-1)$.

 In addition, all the reported simulation results are based on 500 simulations with a nominal significance level of 0.05, and we set the dimension to $100$.

\subsection{Fluctuations of test}

 We illustrated the theoretical asymptotic normality result on simulated data for testing the two-sample problem (\ref{hy1}), and we set precision matrices equal under a null hypothesis, i.e., $\Theta_0^{[1]}=\Theta_0^{[2]}$.

 Letting $G$ be a $p\times p$ symmetric graph matrix with diagonal entries $0$ and $\tilde \alpha$ percent of off-diagonal elements $1$, and $U$ be $p\times p$ matrix with elements i.i.d. generated from the uniformly distribution on the interval $(0,1)$, i.e., $U(0,1)$, we
 denote the elements of the symmetric matrix $\widetilde{\bf \Theta}$ as $\tilde\theta_{ij}$. For $i>j$,
\bqa
\tilde\theta_{ij}=\frac{g_{ij}u_{ij}+g_{ji}u_{ji}}{2}-{\bf 1}_{\{\frac{g_{ij}u_{ij}+g_{ji}u_{ji}}{2}<0.5\}},
\eqa
 where $g_{ij}$ and $u_{ij}$ are the $(i,j)$-entry of $G$ and $U$, respectively, and ${\bf 1}_{\{\cdot\}}$ is the indicator function.
 For $i<j$, we set $\tilde\theta_{ij}=\tilde\theta_{ji}$.
 The diagonal entries of matrix $\widetilde{\bf \Theta}$ are zeros.
 Then, the precision matrix is generated as
\bqa\label{Thetagenerate}
\Theta_0^{[k]}=\widetilde{\bf \Theta}+\left(|\Lambda_{\mathrm{min}}(\widetilde{\bf \Theta})|+0.1\right)\bbI_p.
\eqa

 This shows that the matrix generated is symmetric and positive definite.
 To make the non-zero entries go away from $0$ and to generate a sparse matrix, we subtract $1$ from the non-zero elements.
 In addition, the precision matrix generation procedure shows that $\tilde \alpha$ is a parameter controlling the sparsity.
 When $\tilde \alpha=1$, a dense matrix is generated.
 As is well known, the sparsity of a matrix not only requires a small quantity of non-zero elements, but also a large absolute value of non-zero elements.
 The parameter $\tilde \alpha$ controls sparsity in terms of the number of sparse elements.

 We examined the fluctuation of $\sqrt{n}T_{ij}/\hat{\sigma}_{ij}$ under $(p,n)=(100,200)$ and $(p,n)=(100,400)$ settings for the extremely sparse and dense precision matrix cases, respectively.
 For the extremely sparse precision matrix case, we set the parameter $\tilde\alpha=0.01$, and for dense case we use $\tilde\alpha=1$.

\begin{figure}[htbp]
	\begin{center}
        \includegraphics[width=3.2cm,height=3cm]{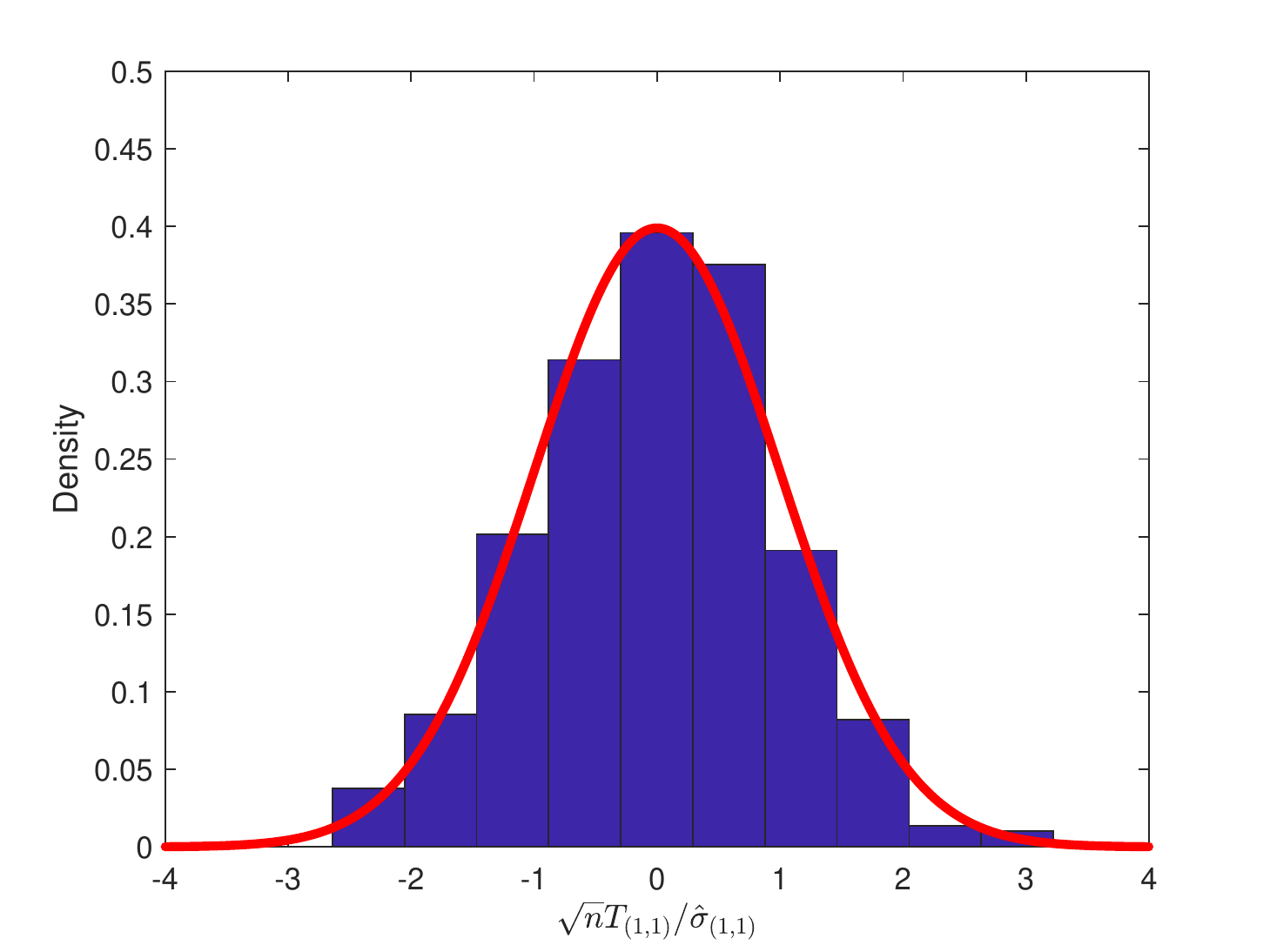}
        \includegraphics[width=3.2cm,height=3cm]{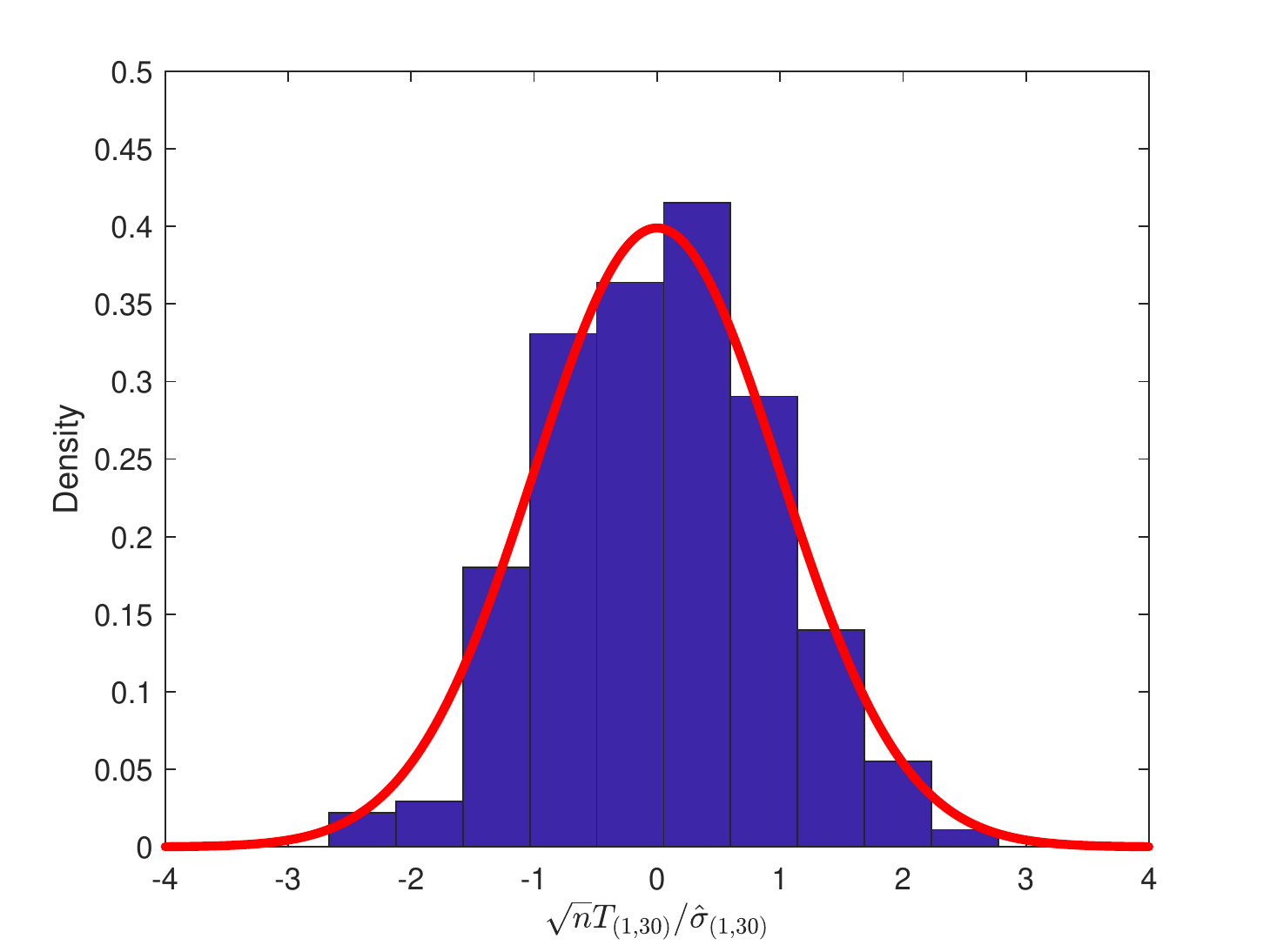}
        \includegraphics[width=3.2cm,height=3cm]{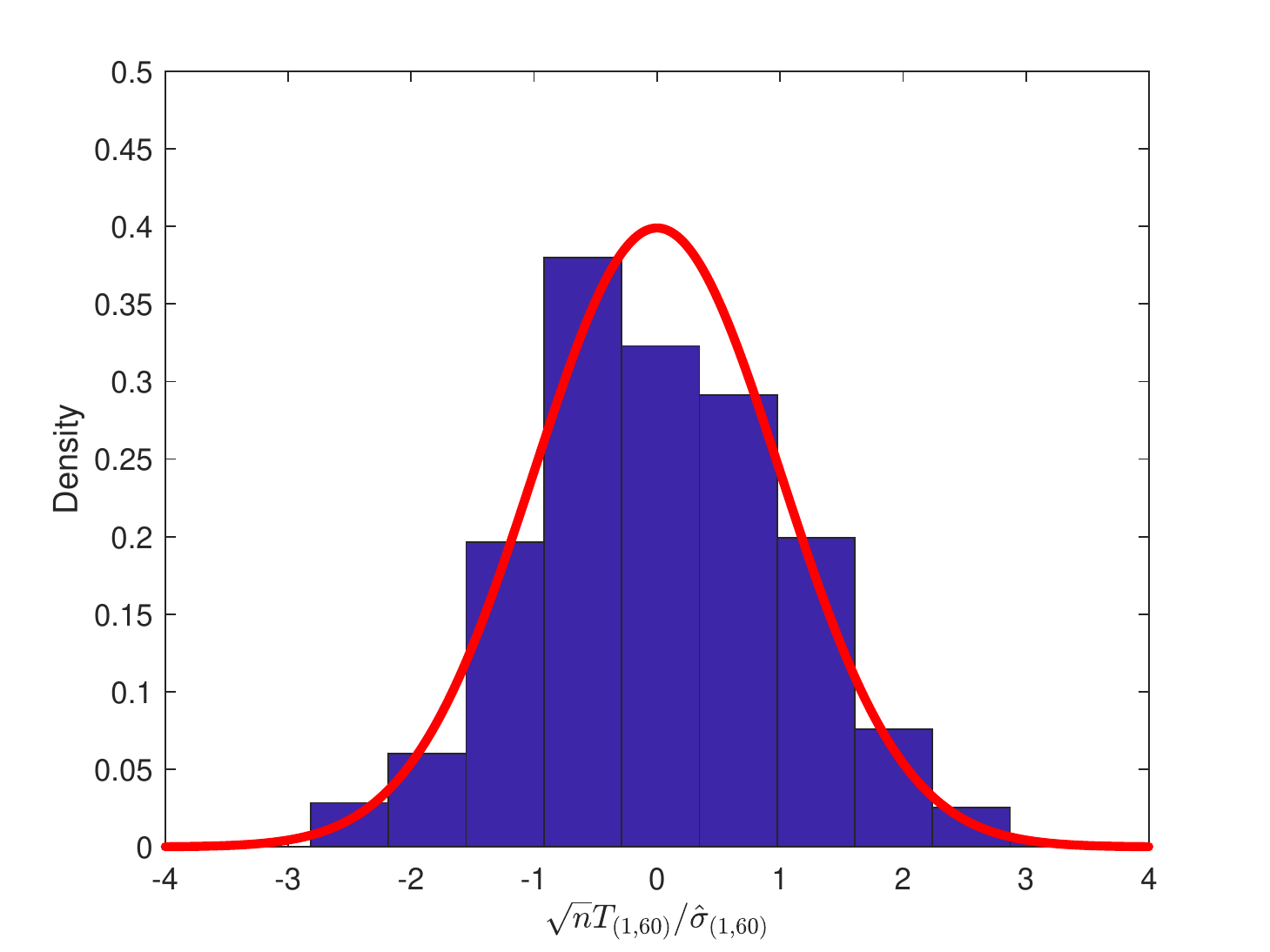}
        \includegraphics[width=3.2cm,height=3cm]{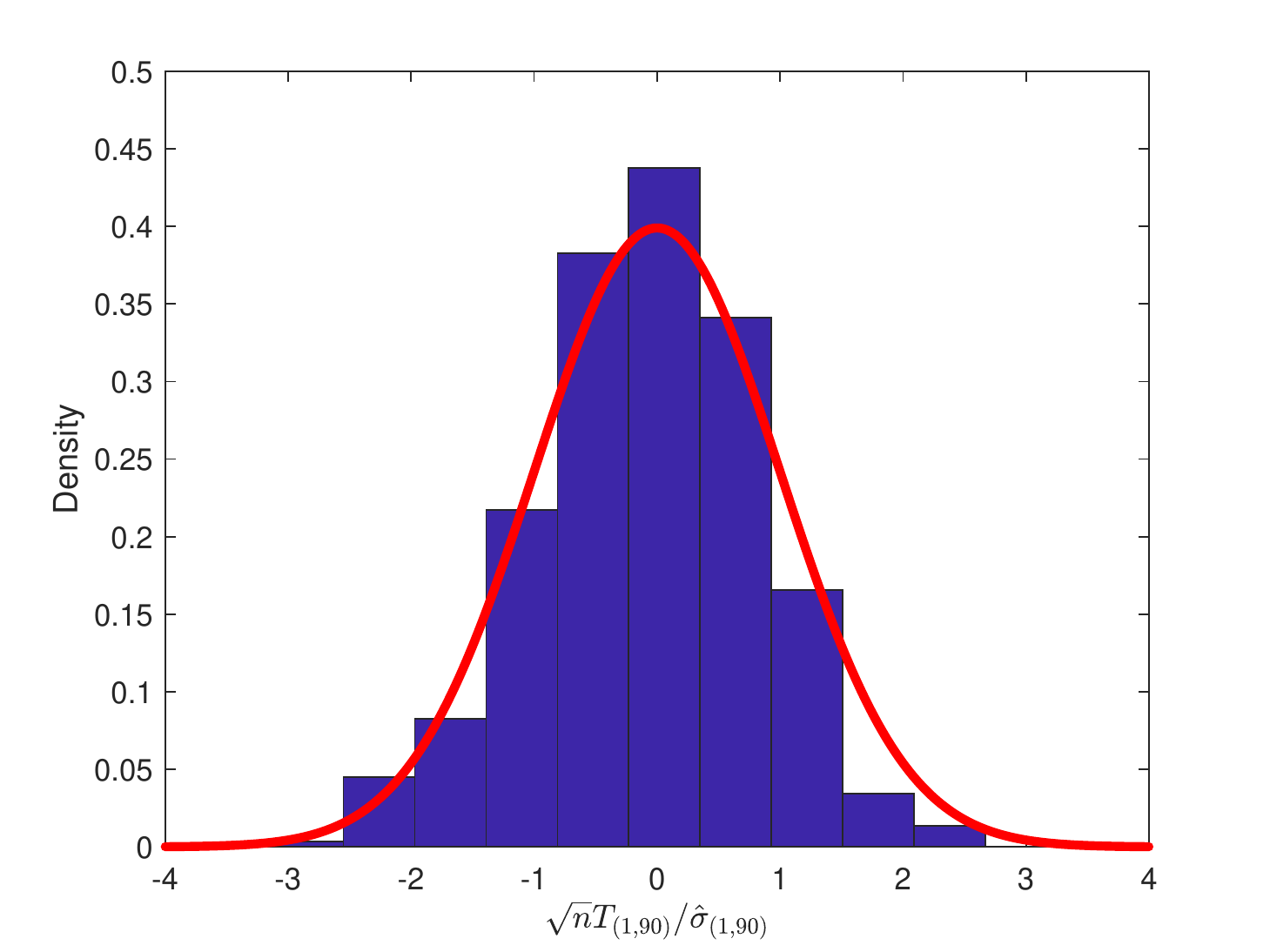}
        \includegraphics[width=3.2cm,height=3cm]{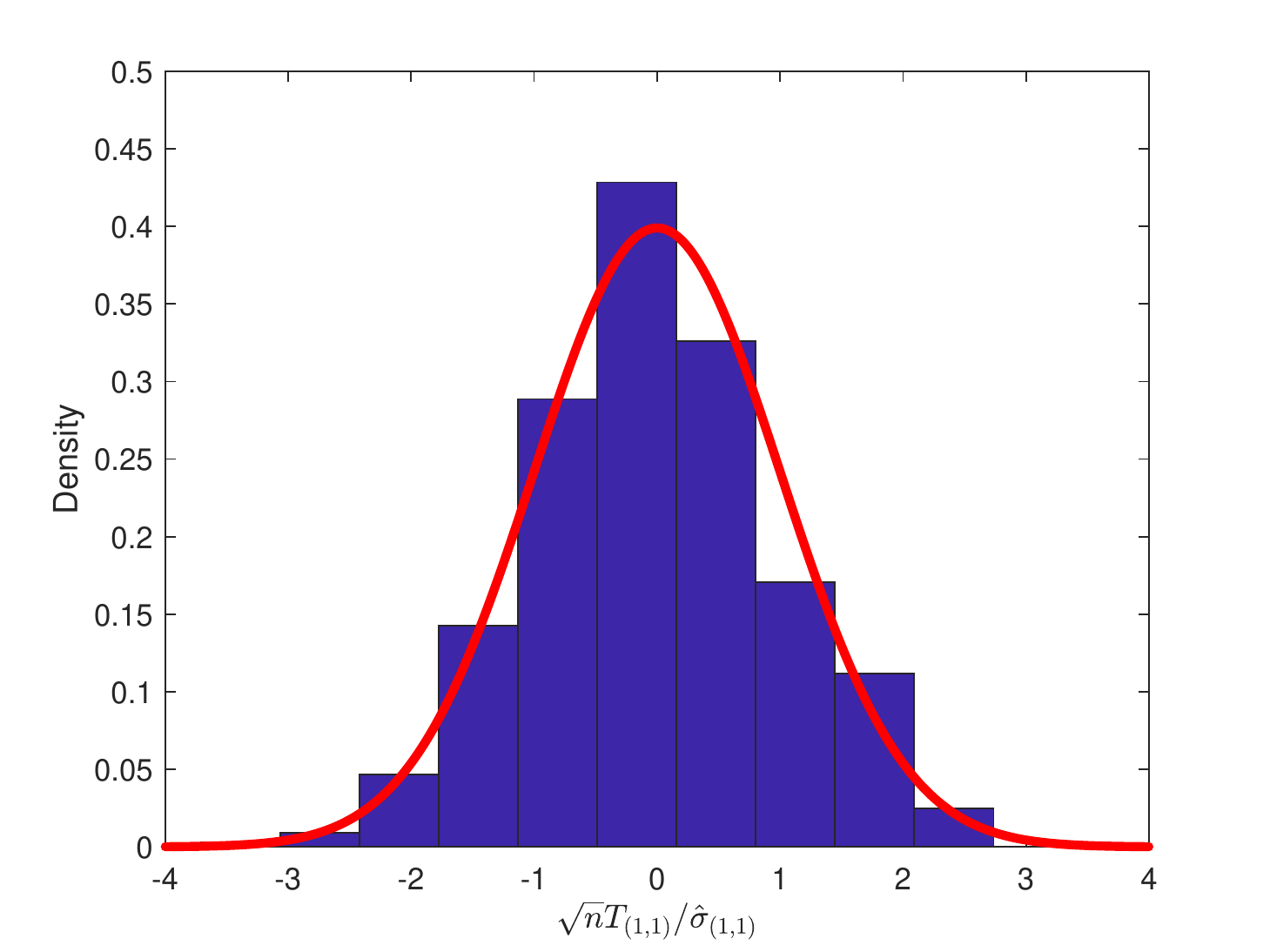}
        \includegraphics[width=3.2cm,height=3cm]{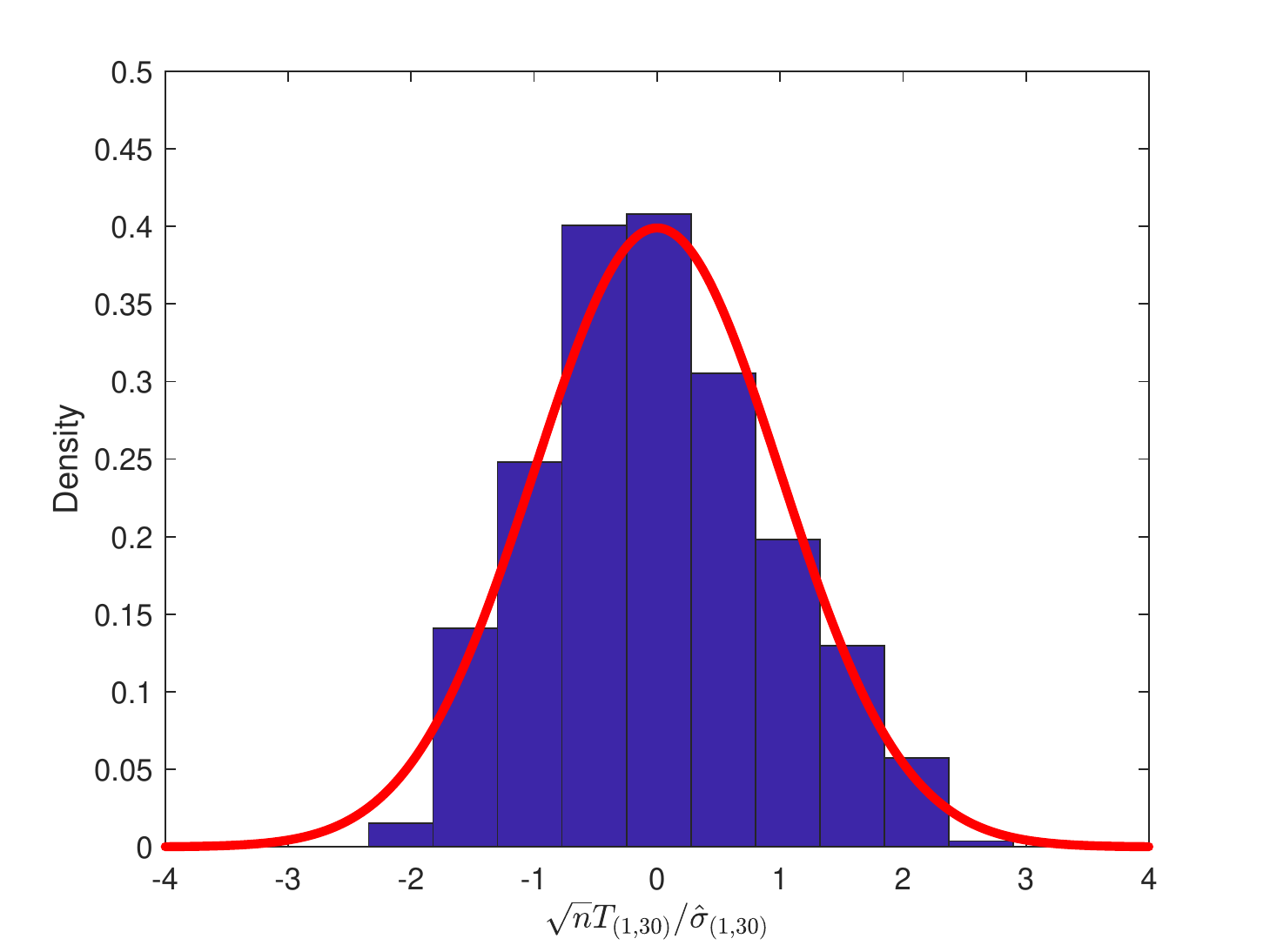}
        \includegraphics[width=3.2cm,height=3cm]{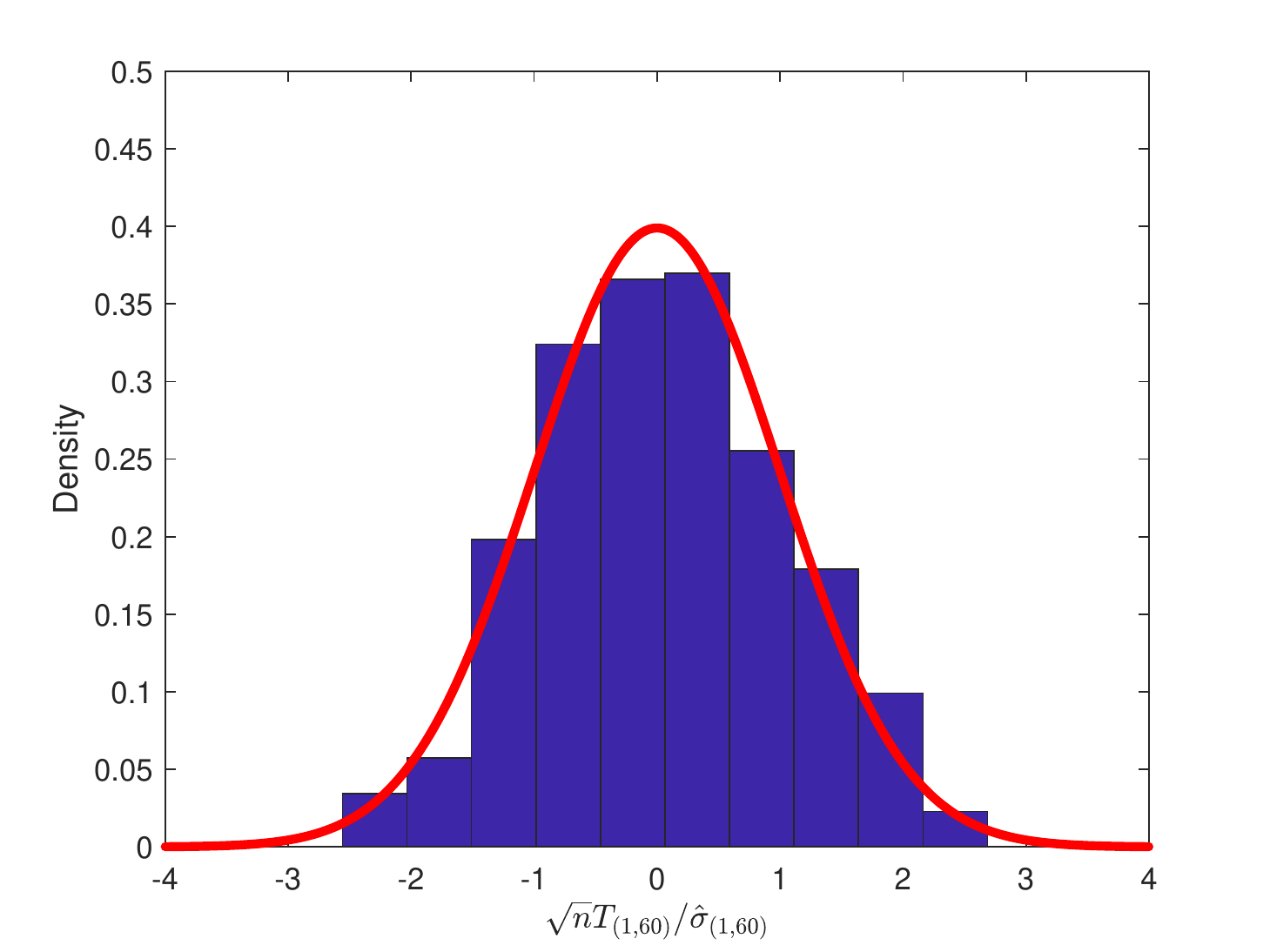}
        \includegraphics[width=3.2cm,height=3cm]{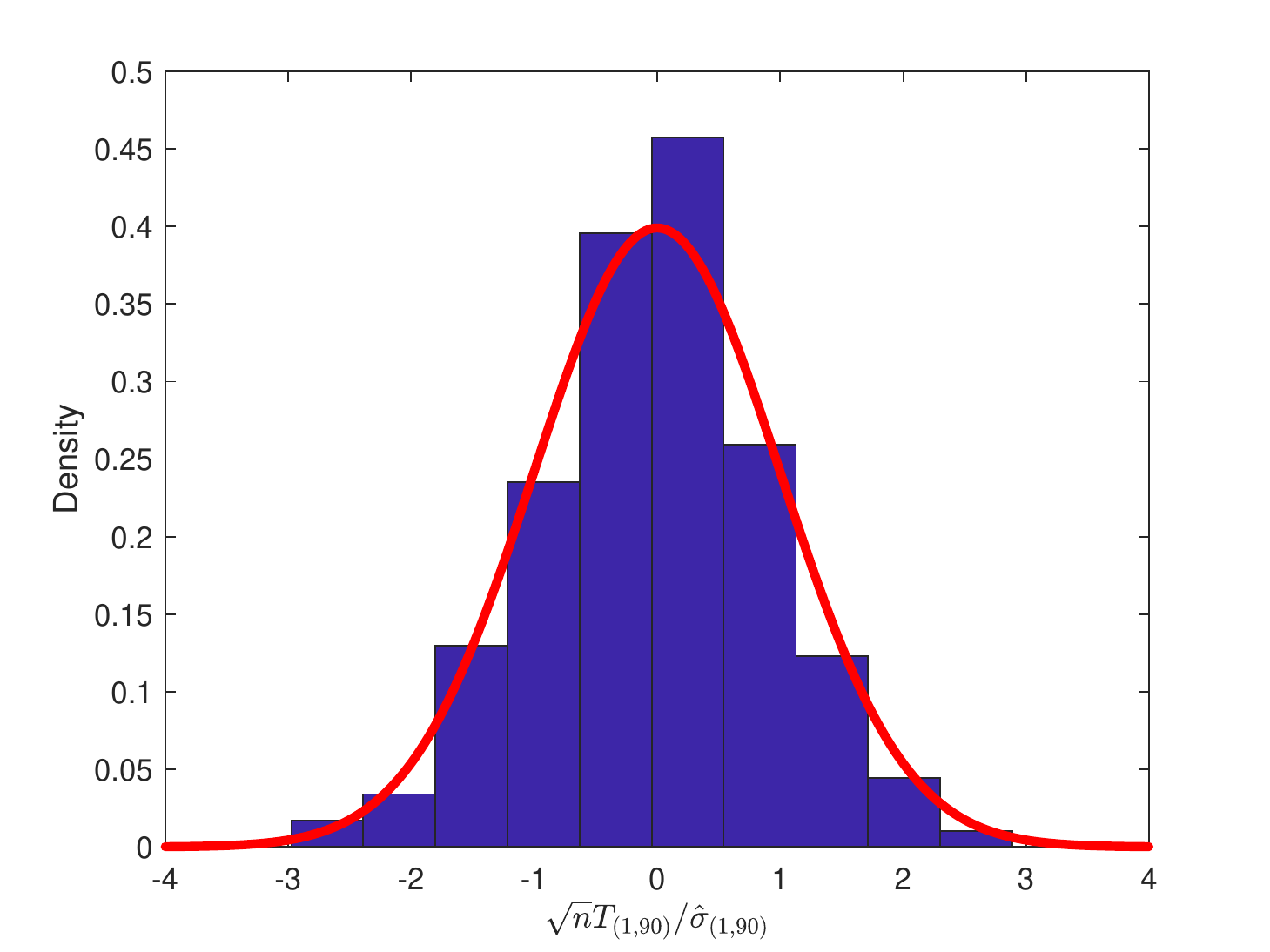}
  \caption{\label{Fluctuation001}Histogram of $\sqrt{n}T_{ij}/\hat{\sigma}_{ij}$ for $\tilde\alpha=0.01$. Here, $T_{(i,j)}=T_{ij}$ and $\hat{\sigma}_{(i,j)}=\hat{\sigma}_{ij}$. The setting is $(p,n)=(100,200)$ with $(i,j)\in\{(1,1), (1,30), (1,60), (1,90)\}$ for four graphs in the first line. The sample size and dimension were set as $(p,n)=(100,400)$ for four graphs in the second line.}
        \end{center}
\end{figure}

\begin{figure}[htbp]
	\begin{center}
        \includegraphics[width=3.2cm,height=3cm]{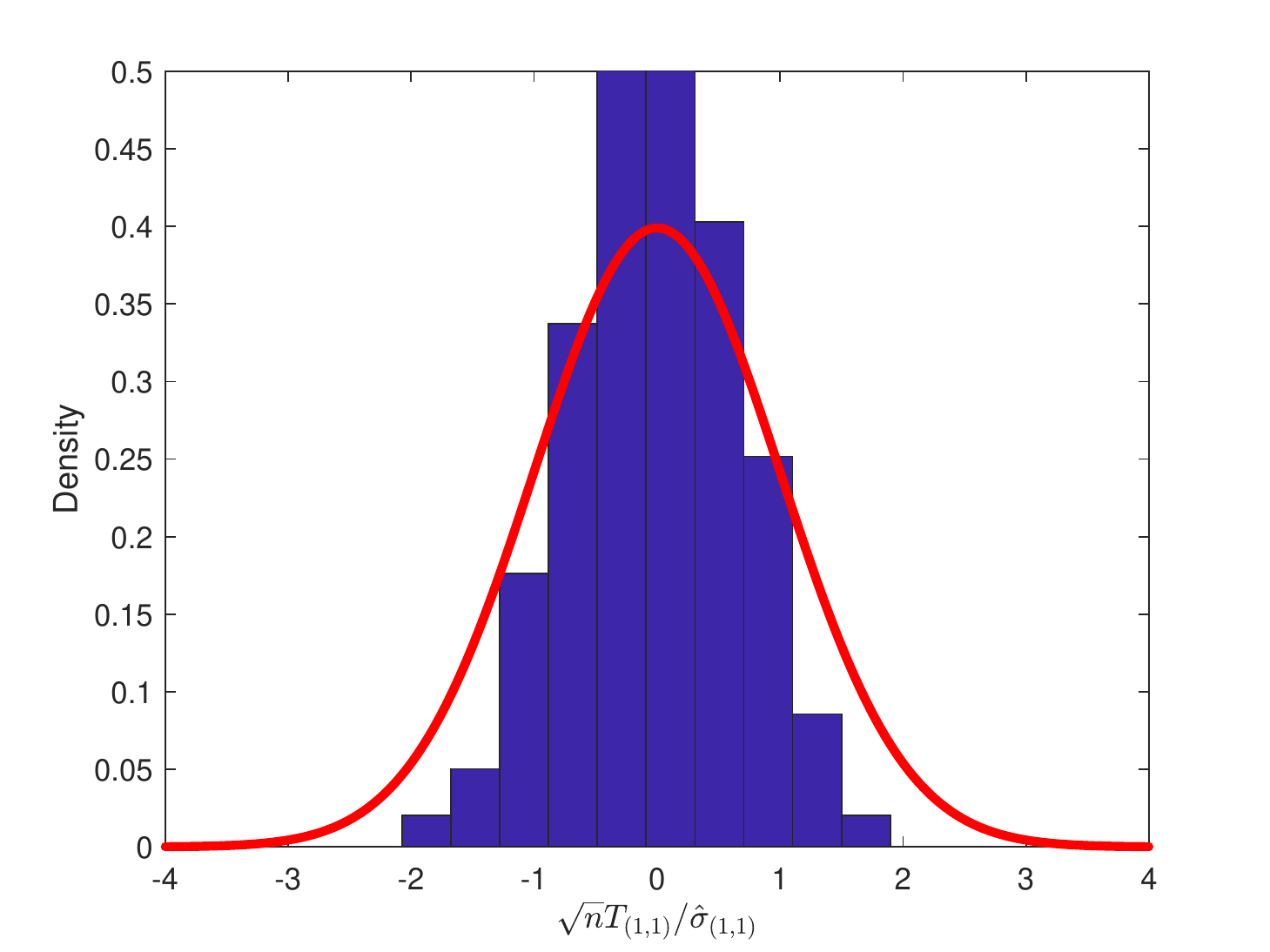}
        \includegraphics[width=3.2cm,height=3cm]{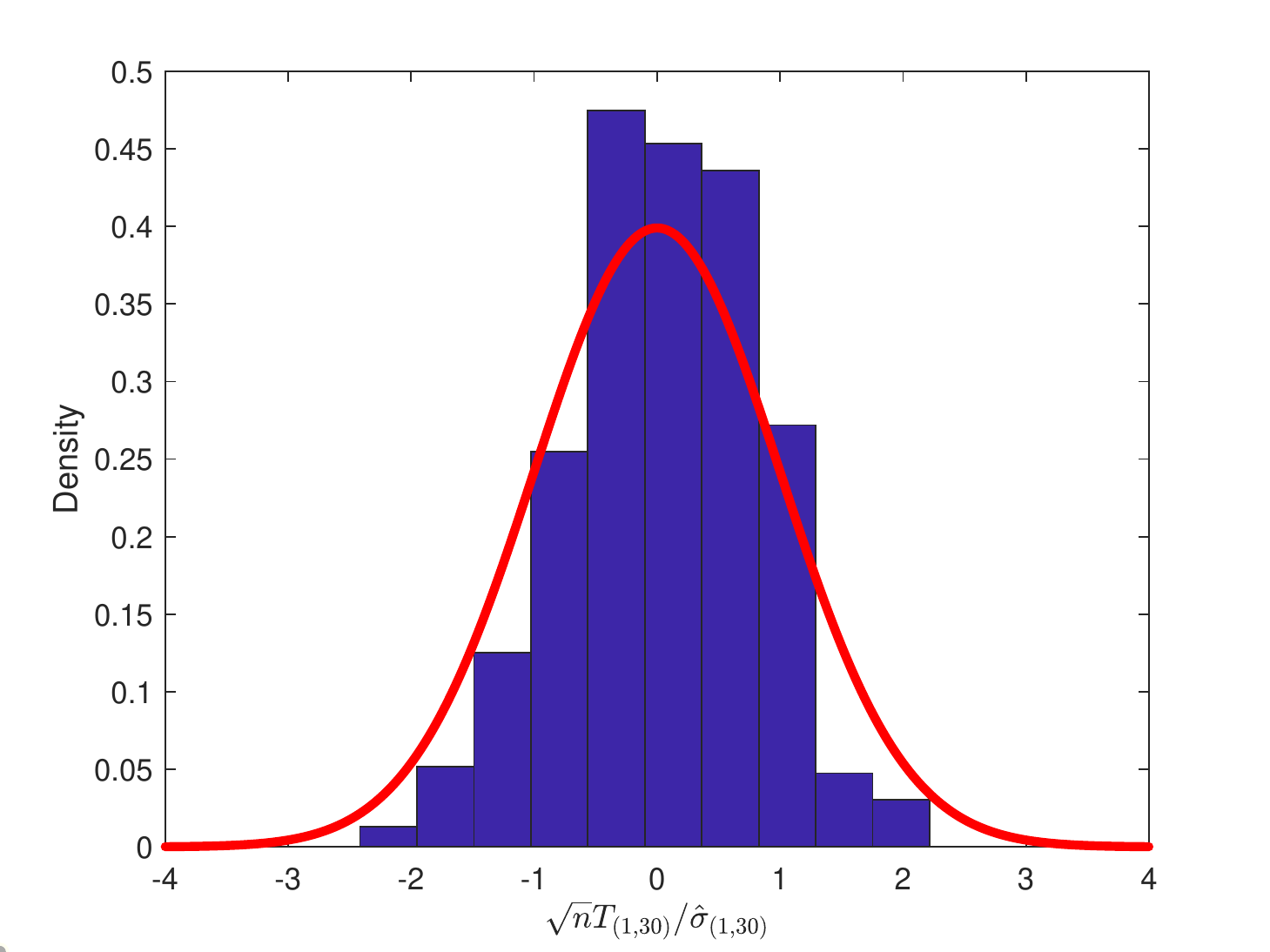}
        \includegraphics[width=3.2cm,height=3cm]{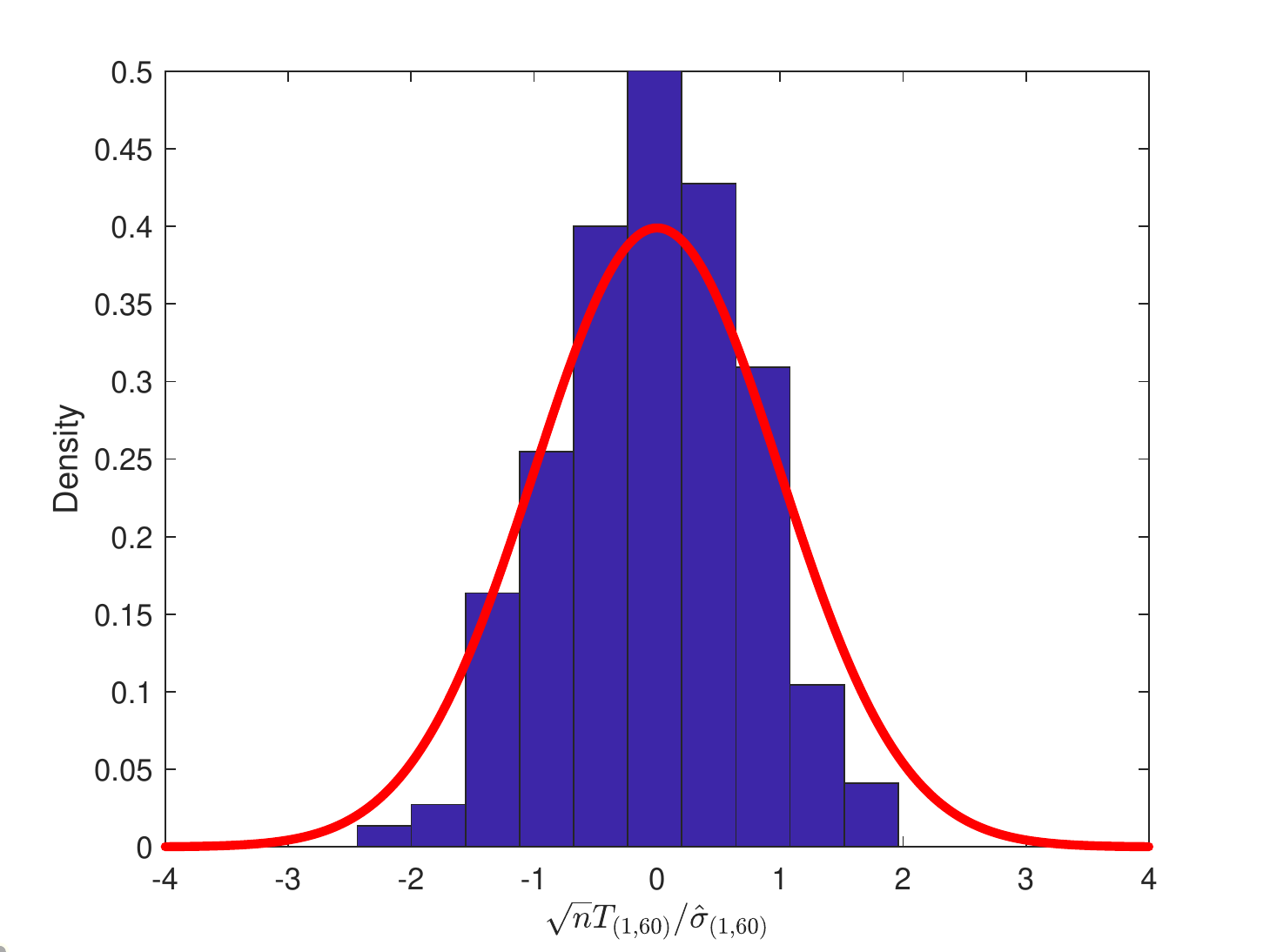}
        \includegraphics[width=3.2cm,height=3cm]{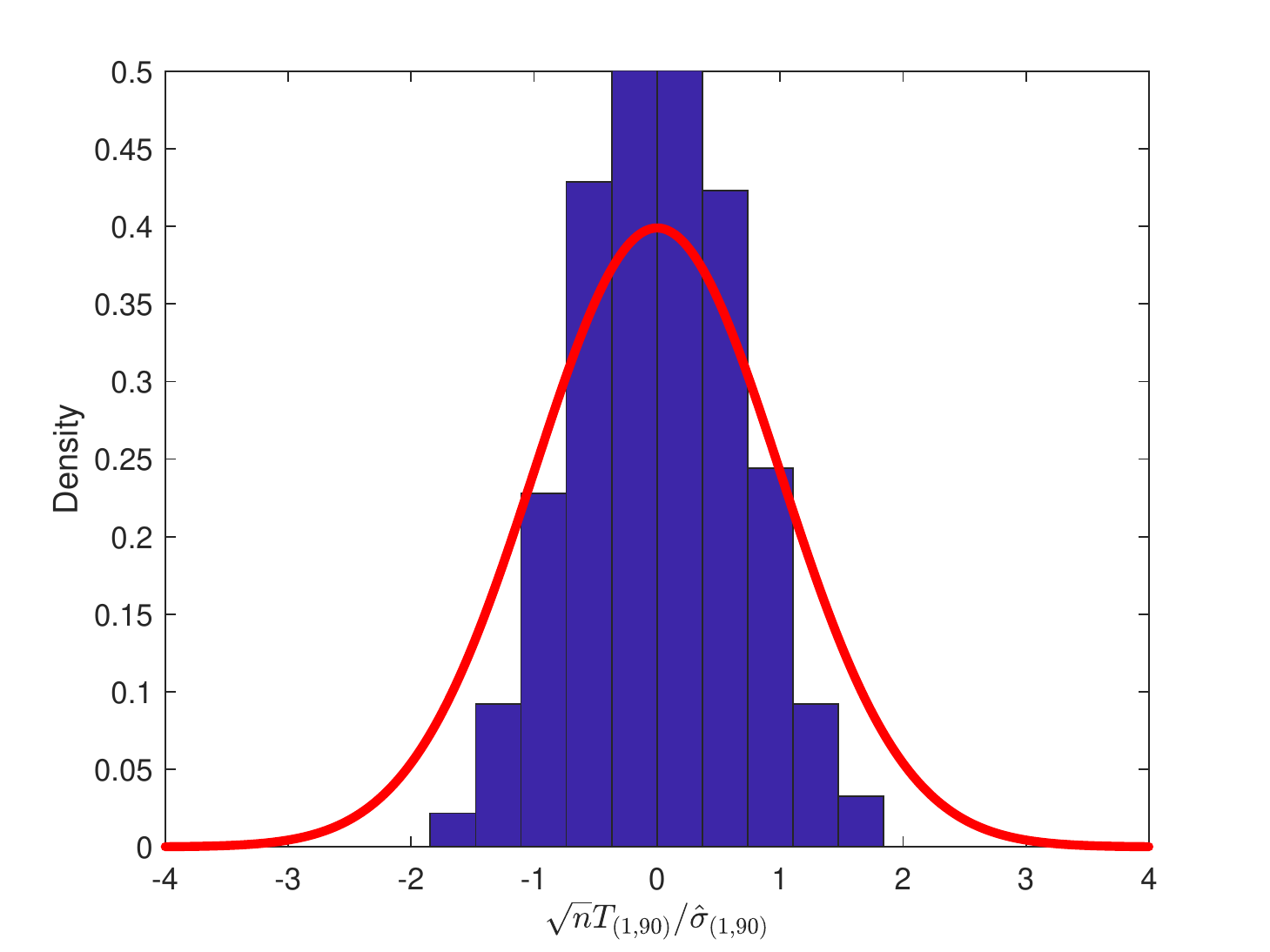}
        \includegraphics[width=3.2cm,height=3cm]{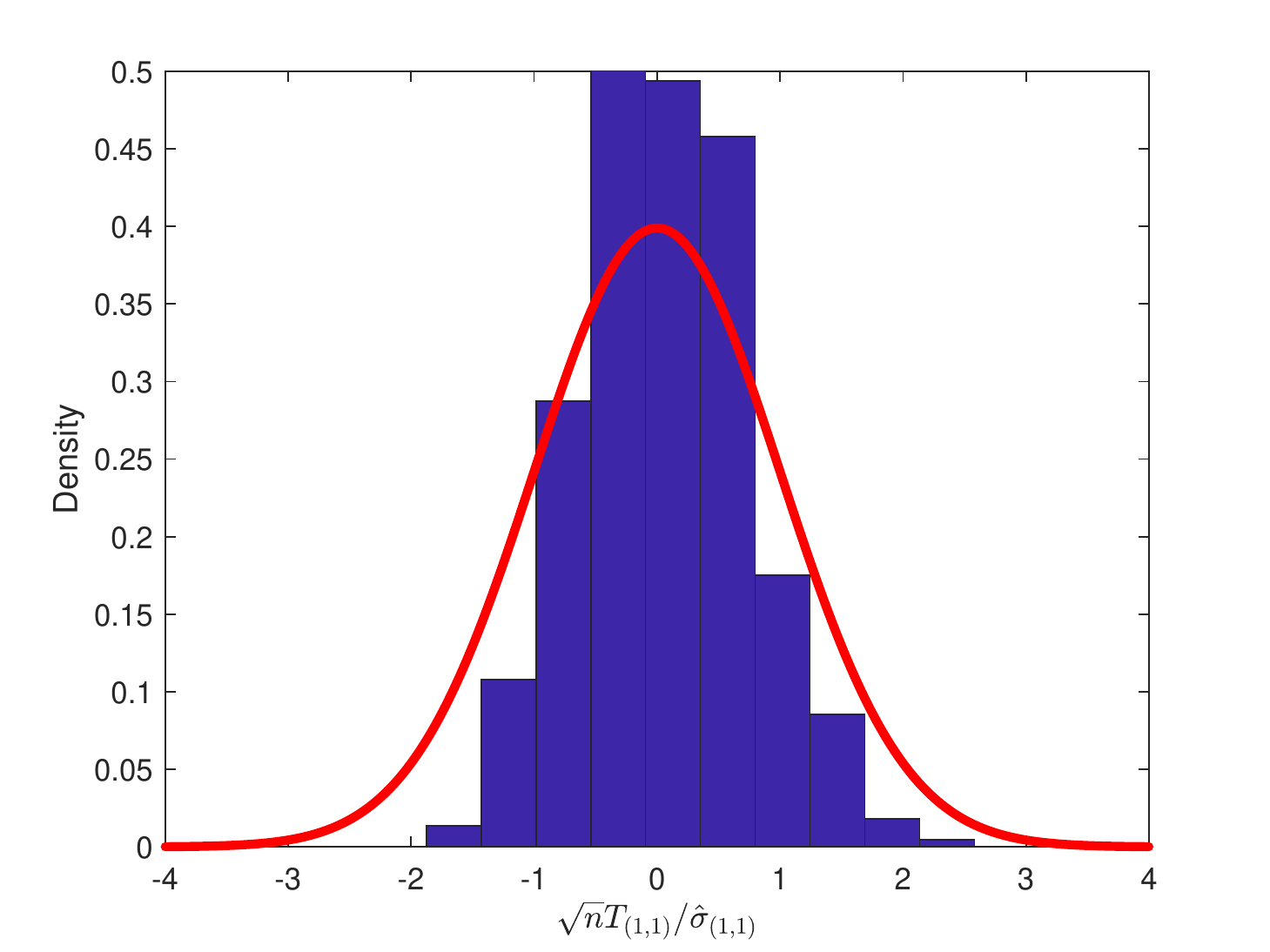}
        \includegraphics[width=3.2cm,height=3cm]{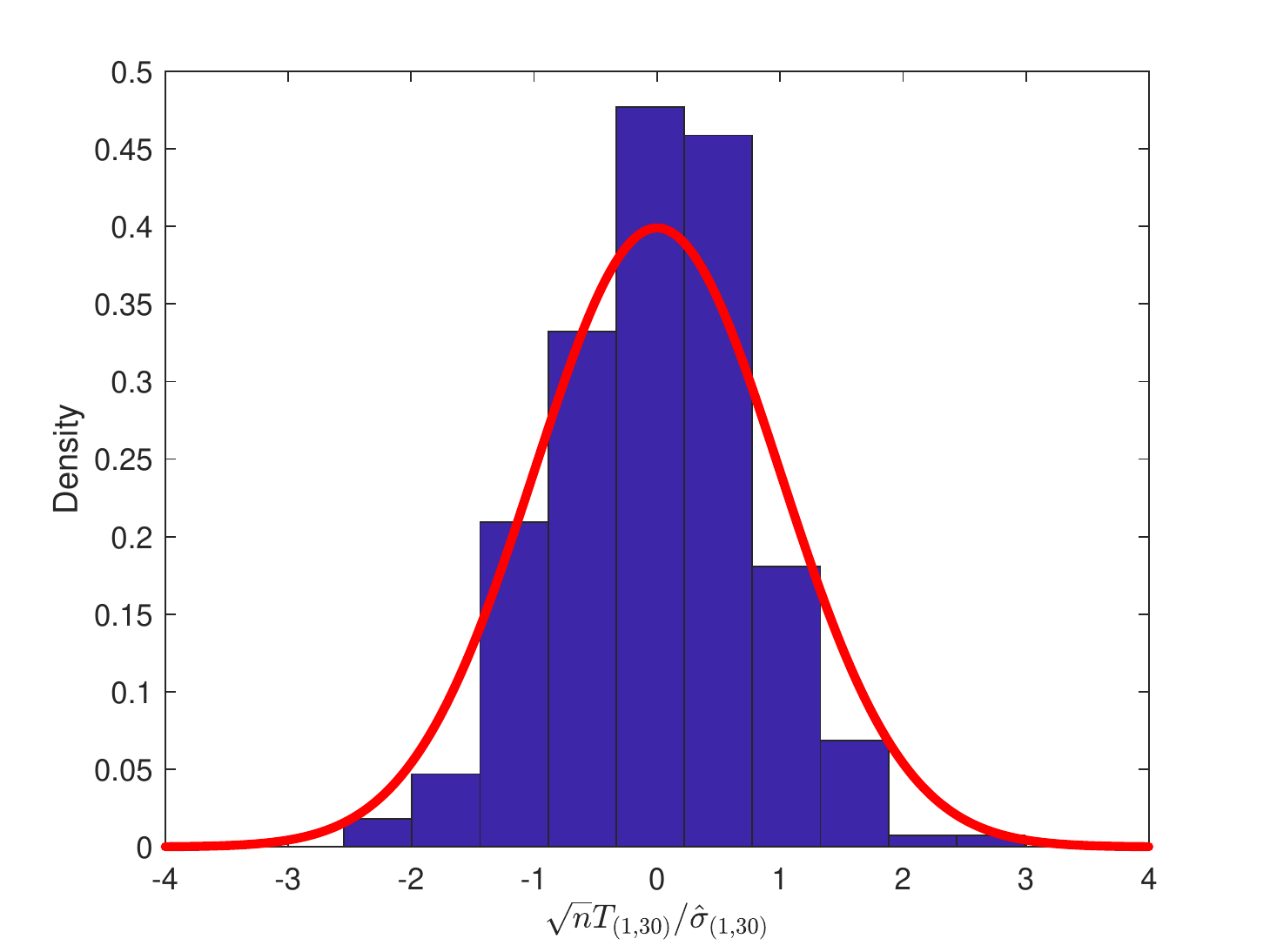}
        \includegraphics[width=3.2cm,height=3cm]{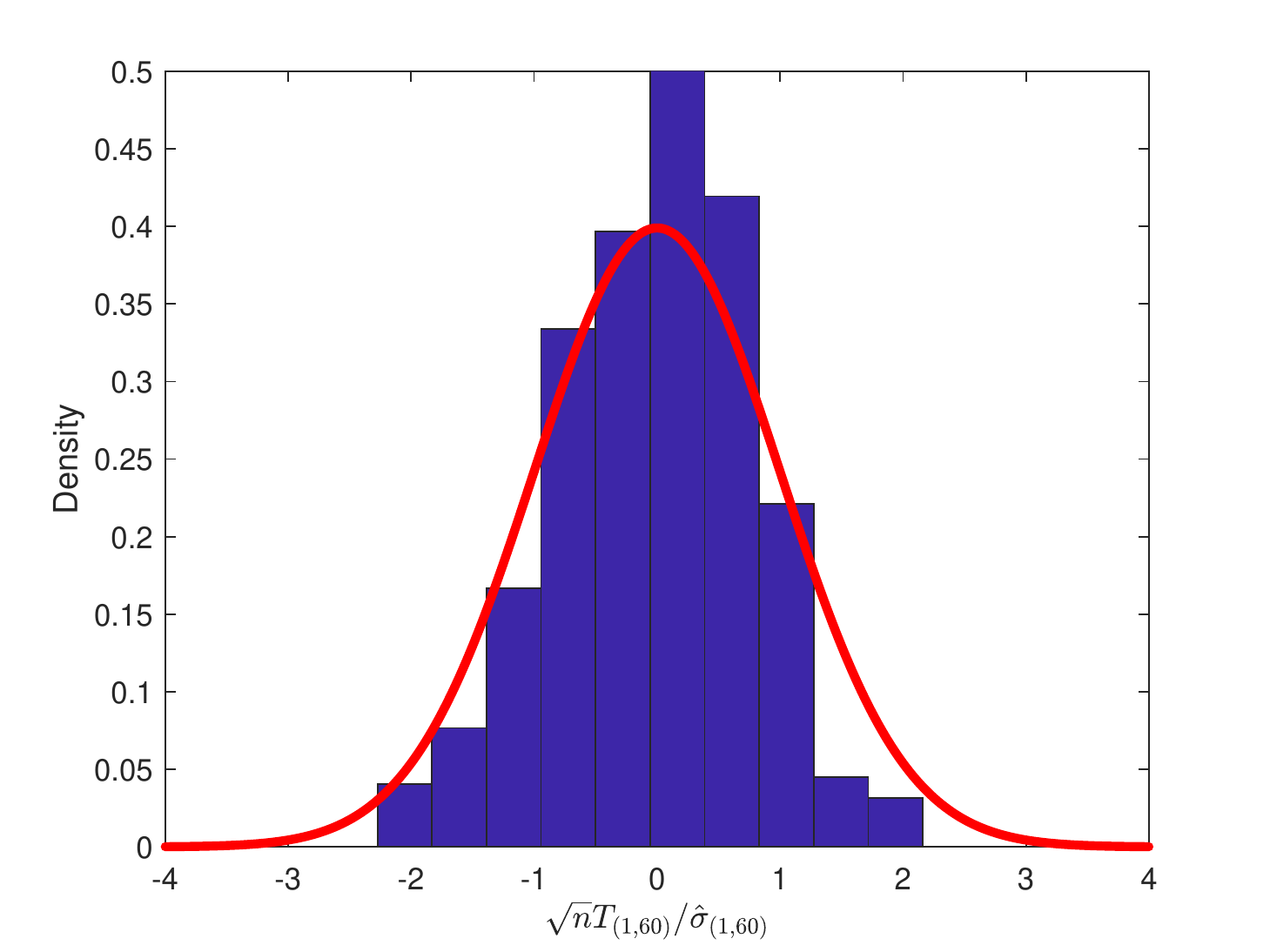}
        \includegraphics[width=3.2cm,height=3cm]{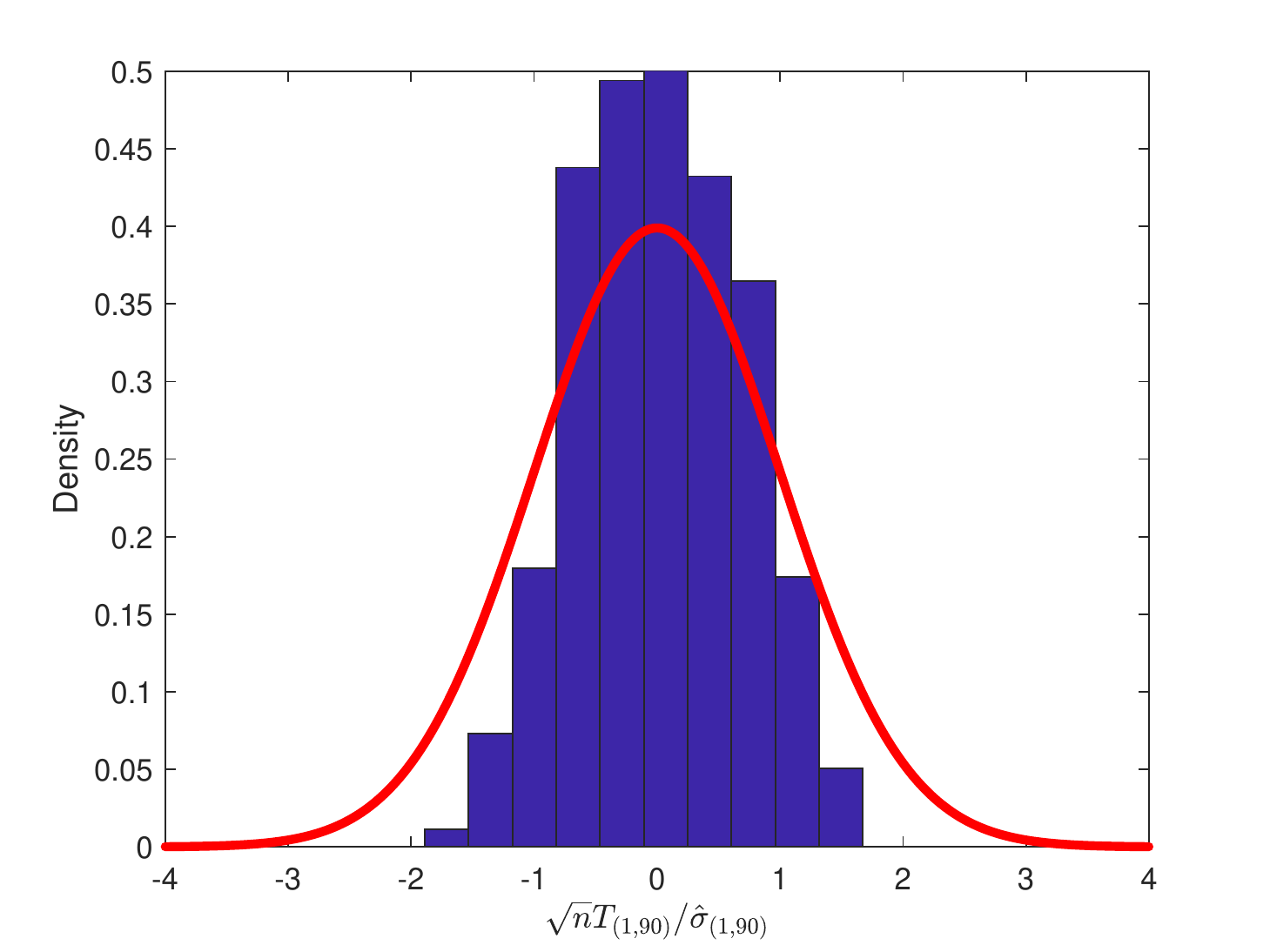}
  \caption{\label{Fluctuation1}Histogram of $\sqrt{n}T_{ij}/\hat{\sigma}_{ij}$ for $\tilde\alpha=1$. Here, $T_{(i,j)}=T_{ij}$ and $\hat{\sigma}_{(i,j)}=\hat{\sigma}_{ij}$. The setting is $(p,n)=(100,200)$ with $(i,j)\in\{(1,1), (1,30), (1,60), (1,90)\}$ for four graphs in the first line. The sample size and dimension were set to $(p,n)=(100,400)$ for four graphs in the second line.}
        \end{center}
\end{figure}

 We simulated the fluctuation for the extremely sparse case as shown in Fig. \ref{Fluctuation001} and the
 dense case in Fig. \ref{Fluctuation1}.
 The index $(i,j)$ in the simulation was intermittently chosen.
 In fact, the CLT provides the method for testing any element of the linear combination of the precision matrix.
 Theoretically, we can test for any index $(i,j)$-entry of ${\bf \Theta}_{0}$ whether the true value is zero or not.

\subsection{Average coverage probabilities}

 We demonstrate the performance of the test method for the $K=2$ situation on testing the hypothesis as follows.
 \begin{itemize}
   \item {\bf Equal Null.} Testing hypothesis (\ref{hy1});
   \item {\bf Linear Null.} Testing the linear null hypothesis $H_0: a_1\Theta^{[1]}_{0ij}+a_2\Theta^{[2]}_{0ij}=0$, i.e., $H_0: \Theta^{[2]}_{0ij}=-\frac{a_1}{a_2}\Theta^{[1]}_{0ij}$. Without loss generation, we chose $-\frac{a_1}{a_2}=0.5$ and $\Theta^{[1]}_{0ij}$ generated from (\ref{Thetagenerate}).
 \end{itemize}

 From the global perspective, we used the average coverage, which is also considered in Jankov$\acute{a}$ and van de Geer \cite{jankova2015confidence}.
 Letting
\bqa\label{confint}
I_{ij}:=\left[T_{ij}-1.96\frac{\sigma_{ij}}{\sqrt{n}}, T_{ij}+1.96\frac{\sigma_{ij}}{\sqrt{n}}\right]
\eqa
 be the $95\%$ asymptotic confidence interval for ${\bf \Theta}_{0ij}$,
 we substitute the estimator $\hat\sigma_{ij}$ for $\sigma_{ij}$ to obtain the empirical version.
 The frequency of the true value being covered by the confidence interval (\ref{confint}) is defined as $\hat\vartheta_{ij}$.
 Then, the average coverage over a set $A$ is denoted
\bqa
Avgcov_{A}=\frac{1}{|A|}\sum_{(i,j)\in A}\hat\vartheta_{ij}.
\eqa
 $S$ denotes the set of non-zero entries of $\Theta^{[1]}_{0ij}$.
 It is easy to check that $S=S_1=S_2$ for the reason that $\Theta^{[1]}_{0ij}$ and $\Theta^{[2]}_{0ij}$ have same structure of sparsity for the Equal Null and Linear Null cases.
 Thus, for the different null hypotheses, we simulated the average coverage over $S$ and its complementary set $S^c$.
 The parameter of sparsity is $\tilde\alpha=0.1, 0.5,$ and $0.9$.

\begin{table}[htbp]
\begin{center}
\setlength{\tabcolsep}{5mm}{
\caption{Estimated average coverage probabilities for $K=2$ situation.}
\label{EmpiricalSizes}
\begin{tabular}{llllll}\hline\noalign{\smallskip}
\multirow{2}*{\rm $\tilde\alpha$}& \multirow{2}*{\rm $n$}& \multicolumn{2}{c}{{\rm Equal Null}} & \multicolumn{2}{c}{{\rm Linear Null}}
\\ \cline{3-6}
& & $S$ & $S^c$ & $S$ & $S^c$
\\ \cline{1-6}
\multirow{2}*{\rm 0.1}&200
 & 0.9886 & 0.9875 & 0.9101 & 0.9824
 \\ \cline{2-6}
&400
& 0.9885 & 0.9867 & 0.8607 & 0.9762
\\ \cline{1-6}
\multirow{2}*{\rm 0.5}&200
& 0.9880 & 0.9878 & 0.9384 & 0.9745
\\ \cline{2-6}
&400
& 0.9870 & 0.9868  & 0.8820 & 0.9647
\\ \cline{1-6}
\multirow{2}*{\rm 0.9}&200
& 0.9901 & 0.9899 & 0.9509 & 0.9751
\\ \cline{2-6}
&400
& 0.9889 & 0.9890 &0.9091  & 0.9639
\\
\noalign{\smallskip}
\hline
\end{tabular}
}
\end{center}
\end{table}

 Partial results in Tab. \ref{EmpiricalSizes} meet our expectation.
 However, we do not deny that the simulations are affected by randomness.
 In addition, the proposed method is based on the combination of estimation and hypothesis testing, which accumulates error.
 The simulation results provide guidance for practice.

\subsection{Multiple FGL case}

 For the multiple FGL case, we examined the fluctuation of the statistic $T_{ij}$ for the $K=3$ situation on testing the hypothesis as follows.
 \begin{itemize}
   \item {\bf Three-sample Linear Null.} Testing hypothesis $H_0: \Theta^{[3]}_{0ij}=-\frac{a_1}{a_3}\Theta^{[1]}_{0ij}-\frac{a_2}{a_3}\Theta^{[2]}_{0ij}$,
 where $-\frac{a_1}{a_3}=0.6$ and $-\frac{a_2}{a_3}=0.9$ are both generated from $U(0,1)$.
 $\Theta^{[1]}_{0ij}$ and $\Theta^{[2]}_{0ij}$ are both generated from (\ref{Thetagenerate}) with parameters $0.01$ and $0.1$, respectively.
 \end{itemize}

\begin{figure}[htbp]
	\begin{center}
        \includegraphics[width=3.2cm,height=3cm]{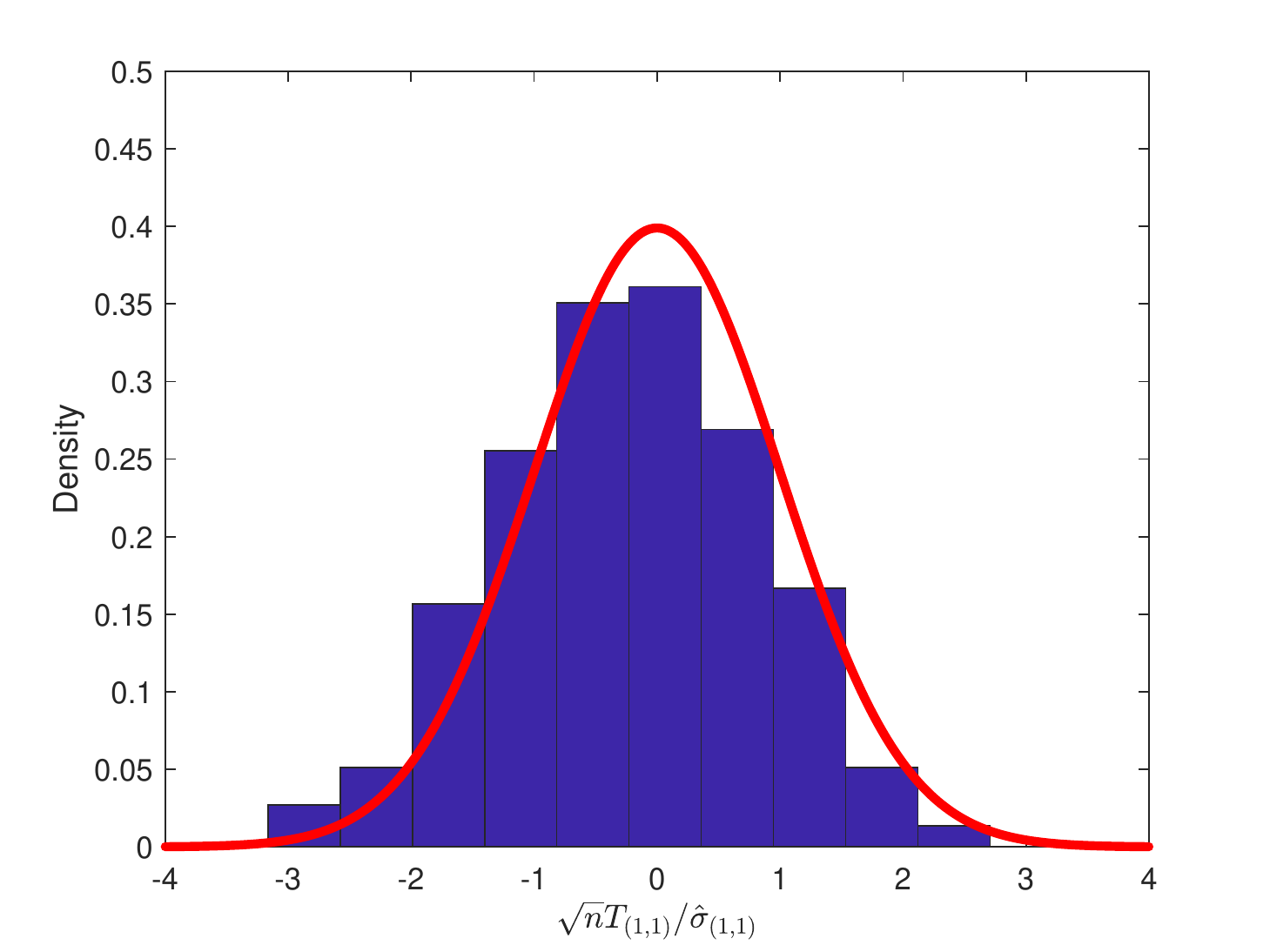}
        \includegraphics[width=3.2cm,height=3cm]{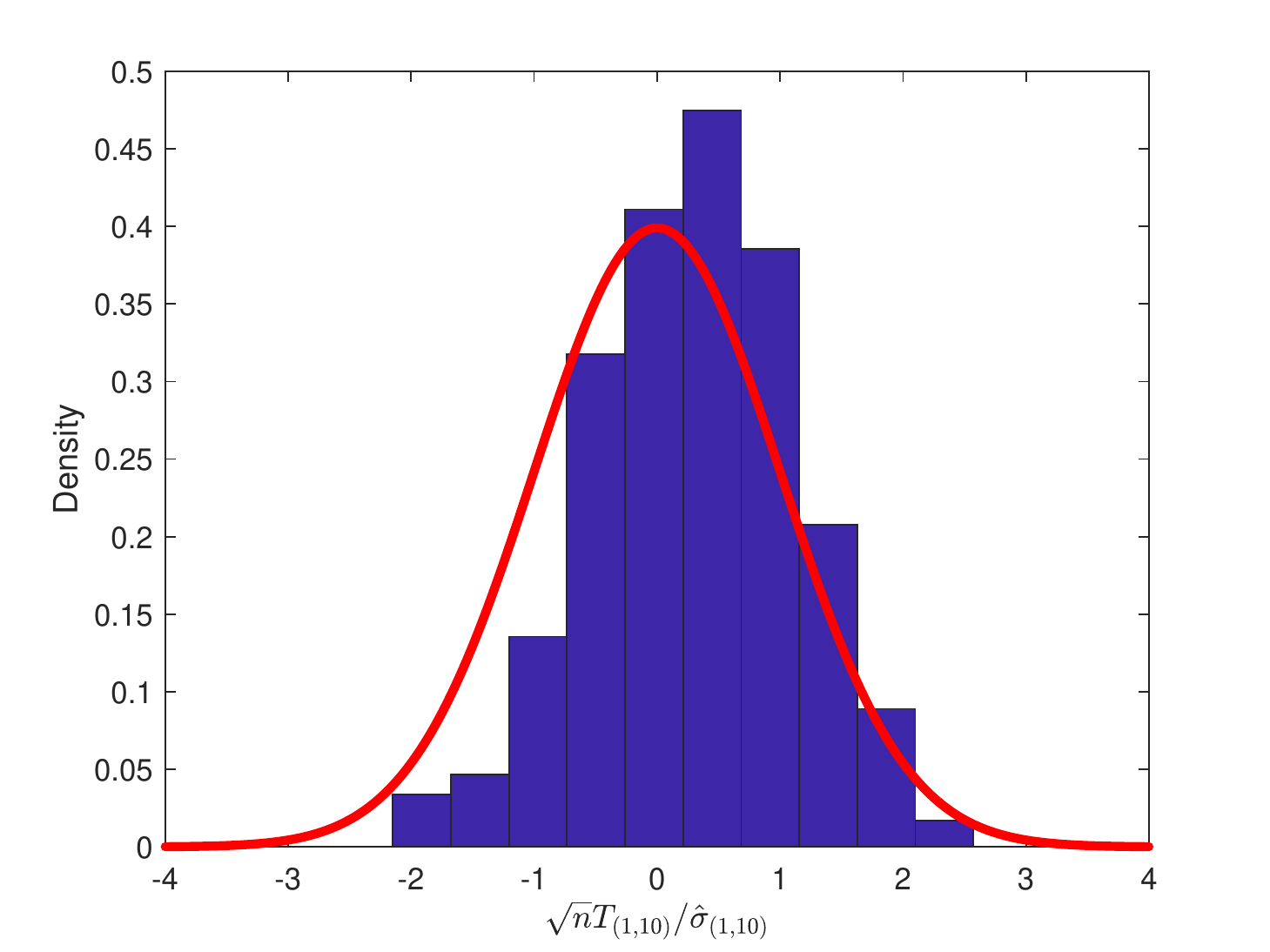}
        \includegraphics[width=3.2cm,height=3cm]{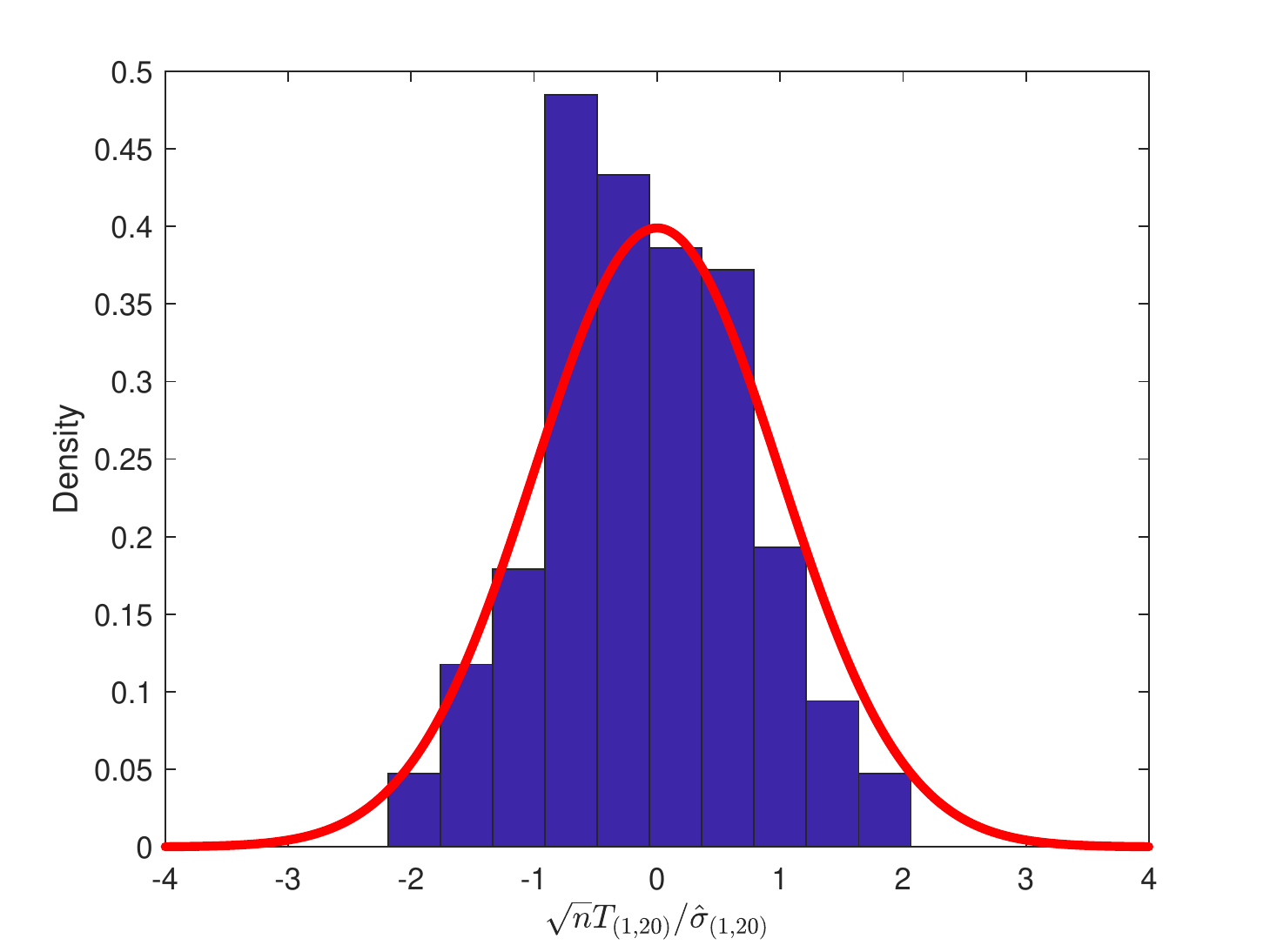}
        \includegraphics[width=3.2cm,height=3cm]{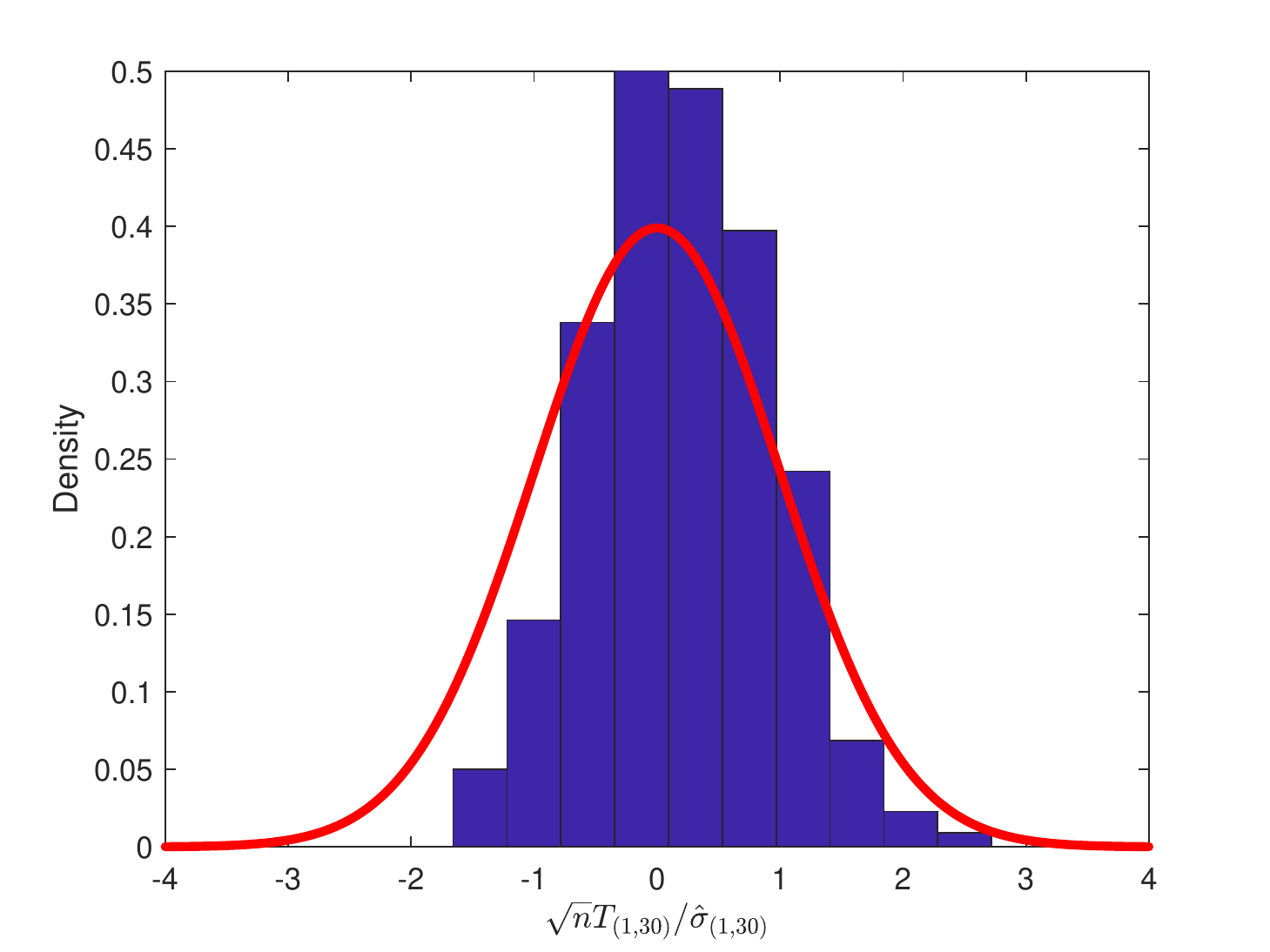}
        \includegraphics[width=3.2cm,height=3cm]{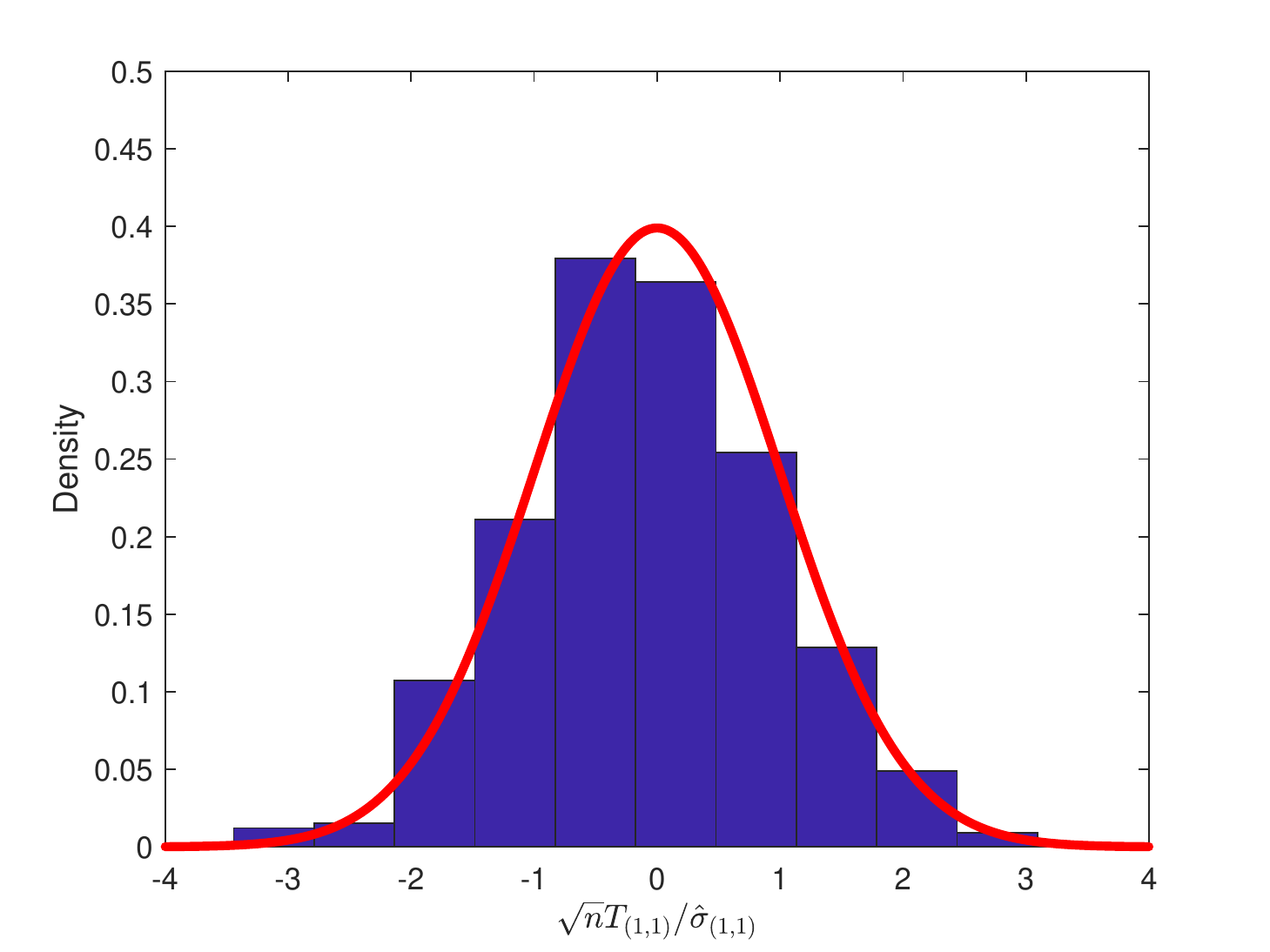}
        \includegraphics[width=3.2cm,height=3cm]{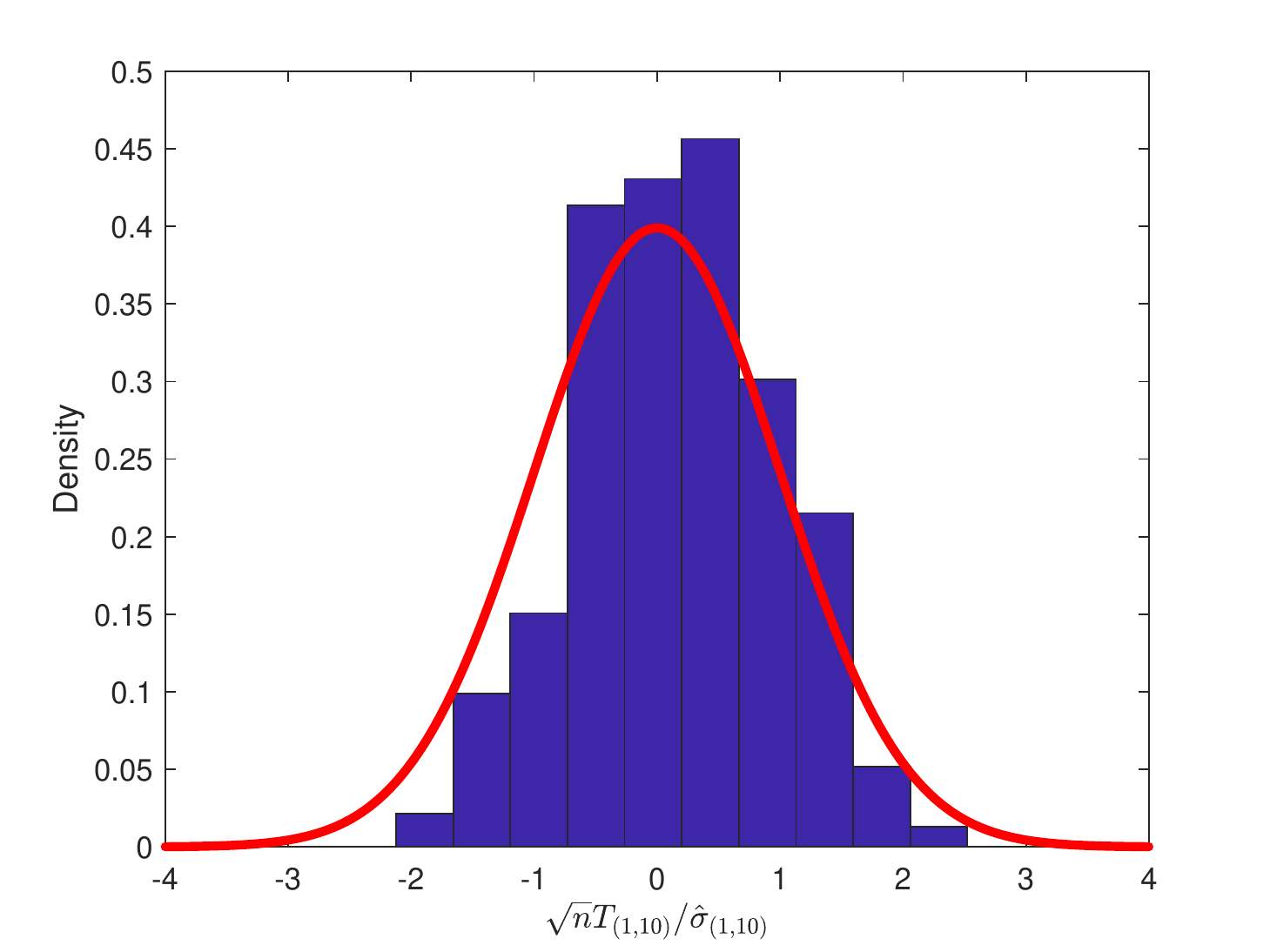}
        \includegraphics[width=3.2cm,height=3cm]{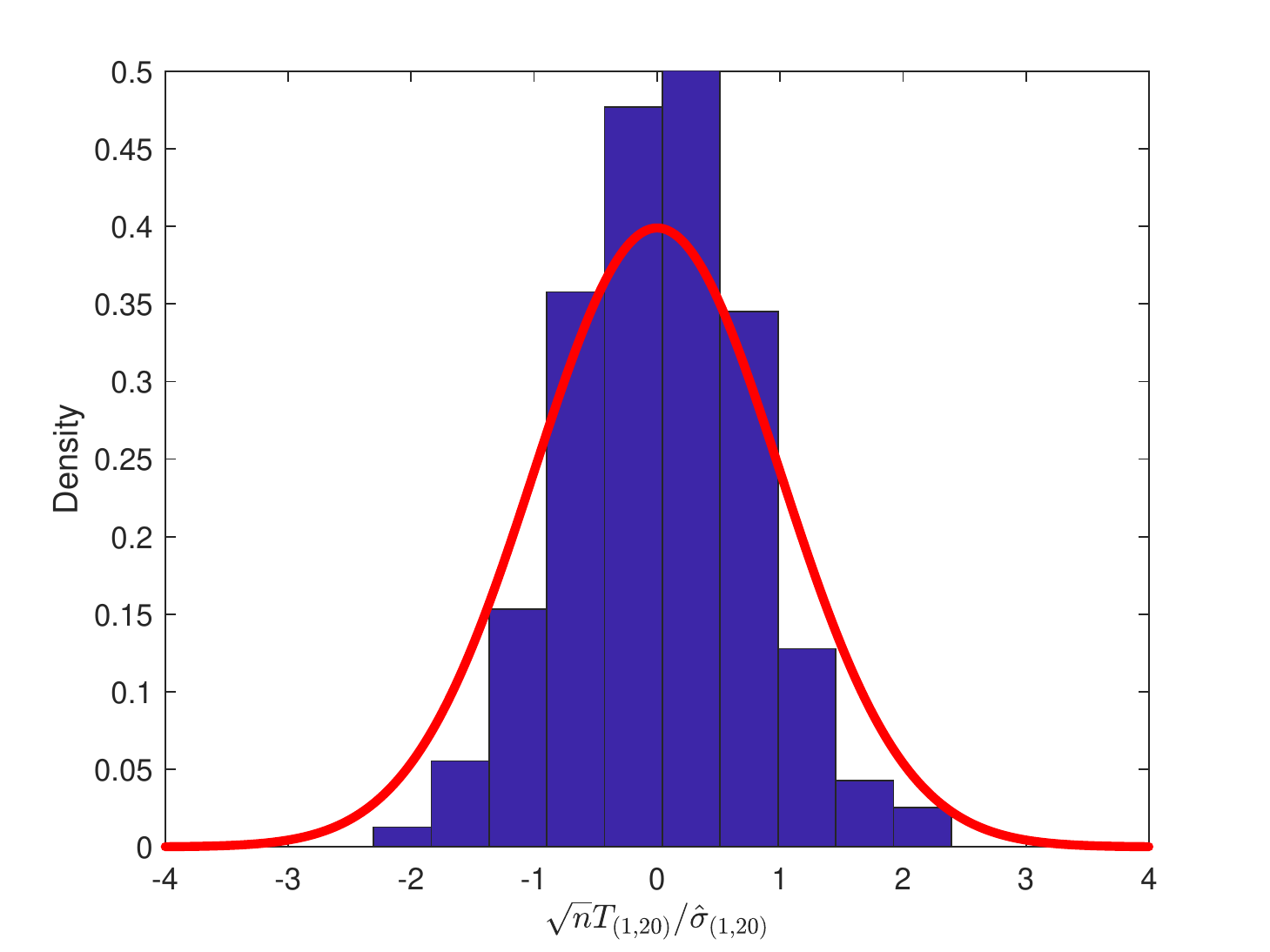}
        \includegraphics[width=3.2cm,height=3cm]{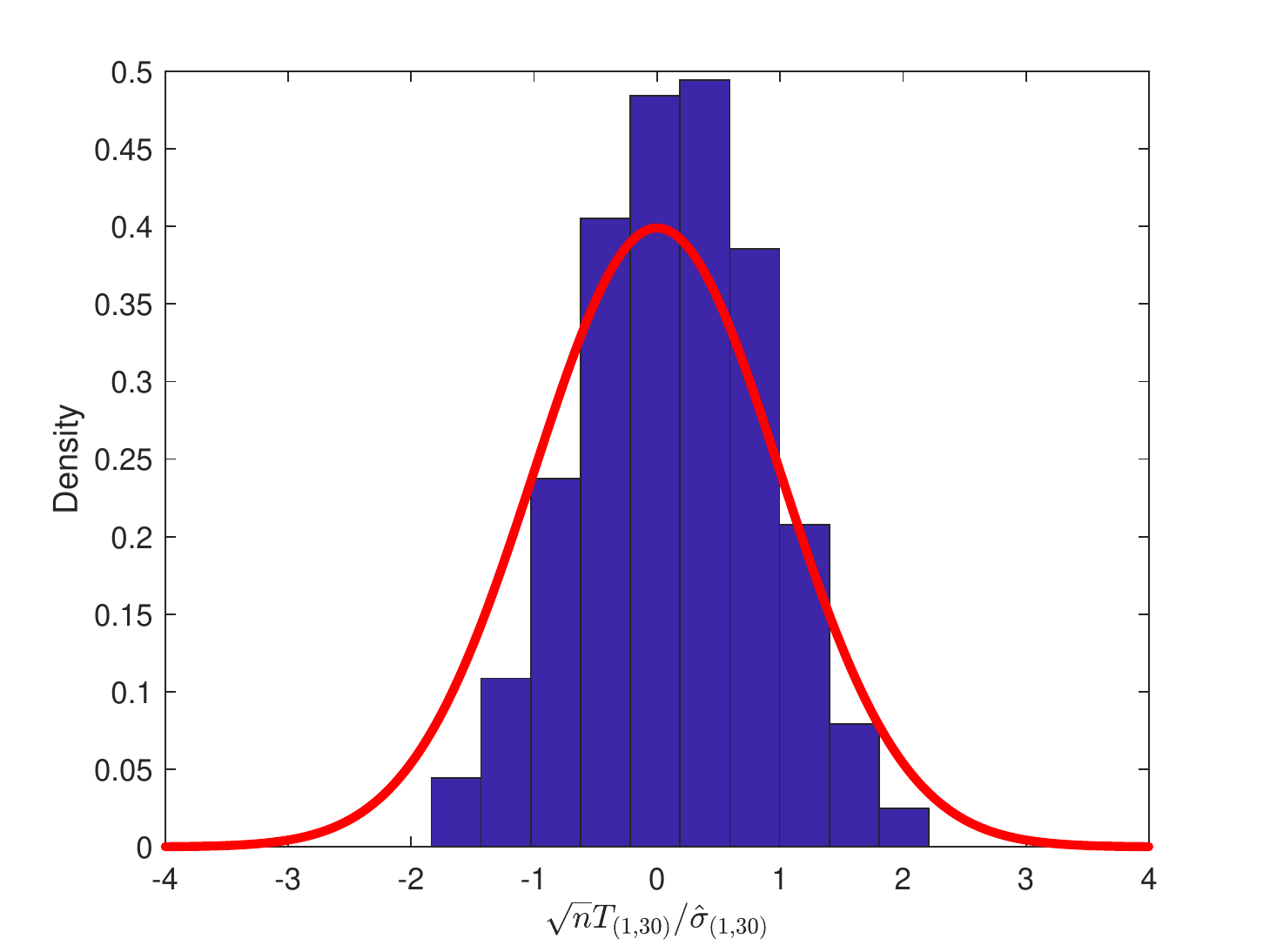}
  \caption{\label{Fluc}Histogram of $\sqrt{n}T_{ij}/\hat{\sigma}_{ij}$ for $\tilde\alpha=1$. Here, $T_{(i,j)}=T_{ij}$ and $\hat{\sigma}_{(i,j)}=\hat{\sigma}_{ij}$. The setting is $(p,n)=(100,200)$ with $(i,j)\in\{(1,1), (1,10), (1,20), (1,30)\}$ for four graphs in the first line. The sample size and dimension were set to $(p,n)=(100,400)$ for four graphs in the second line.}
        \end{center}
\end{figure}

 We set $-\frac{a_1}{a_3}$ and $-\frac{a_2}{a_3}$ to positive numbers, since the setting of hypothesis testing should guarantee that $\{\Theta^{[k]}_{0ij}\}_{k=1}^{3}$ are symmetric positive-definite matrices.
 Besides, for Three-sample Linear Null, $S$ denotes the set of non-zero entries of $\Theta^{[1]}_{0ij}+\frac{a_2}{a_3}\Theta^{[2]}_{0ij}$.
 The dimension and sample size are $(p,n)=(100,200)$ and $(p,n)=(100,400)$, respectively.
 Histograms of the proposed statistic $T_{ij}$ at the
\bqn
(i,j)\in\{(1,1),(1,10),(1,20),(1,30)\}
\eqn
 locations of the precision matrix are presented in Fig. \ref{Fluc}.

\section*{Acknowledgments}

%We thank LetPub (www.letpub.com) for its linguistic assistance during the preparation of this manuscript.

Q.Y. Zhang was partially supported by NSFC 12201430, 11971097, and Capital University of Economics and Business: The Fundamental Research Funds for Beijing Universities XRZ2021044.
Z.D. Bai was supported by NSFC 12171198, 12271536 and STDFJ 20210101147JC.
H. Yang was supported by NSSF China, Grant 22FGLB056.

\section*{References}

\bibliographystyle{elsarticle-harv}
\bibliography{FusedGGM}

\section*{Appendix}

\begin{appendices}

\section{Proof of Theorem}

\subsection{Proof of Theorem \ref{Orcalebounds}}
 To prove Theorem \ref{Orcalebounds}, we need a lemma of Jankov{\'a} and Van de Geer \cite{jankova2018inference}, which is present as follow.

\begin{lem}\label{lembasic}
Let $f(\Delta):=tr(\Delta\Sigma_{0})-[\log\det(\Delta+\Theta_{0})-\log\det(\Theta_{0})]$. Assume that $1/L\leq \lambda_{min}(\Theta_0)\leq \lambda_{max}(\Theta_0)\leq L$ for some constant $L\geq 1$. Then for all $\Delta$ such that $||\Delta||_{F} \leq 1/(2L)$, $f(\Delta)$ is well defined and
\begin{equation}
f(\Delta) \geq \frac{1}{2(L+1/(2L))^2}||\Delta||_{F}^2.
\end{equation}
\end{lem}

 To simplify the notation, we substitute $\widehat{\Sigma}_k$, $\Sigma_{0k}$, $\widehat{\Theta}_k$, $\Theta_{0k}$ for $\widehat{\Sigma}^{[k]}$, $\Sigma_{0}^{[k]}$, $\widehat{\Theta}^{[k]}$, $\Theta_{0}^{[k]}$ respectively.

\begin{proof}

 Note that $\widehat{\Theta}_k$ is the minimum value of the fused graphical Lasso for $k=1,2$.
 Let $\widetilde{\Theta}_k=\alpha_k\widehat{\Theta}_k+(1-\alpha_k)\Theta_{0k}$, and $\alpha_k=\frac{M}{M+||\widehat{\Theta}_k-\Theta_{0k}||_{F}}$. According to the definitions of $\widetilde{\Theta}_k$, and the convexity of loss function $$F_n(\Theta_1,\Theta_2)=tr(\Theta_1\widehat{\Sigma}_1)-\log\det(\Theta_1)
 +tr(\Theta_2\widehat{\Sigma}_2)-\log\det(\Theta_2)+\lambda||\Theta_1^{-}||_1+
 \lambda||\Theta_2^{-}||_1+\rho||\Theta_1^{-}-\Theta_2^{-}||_1,$$ we obtain
\bqa
F_n(\widetilde\Theta_1,\widetilde\Theta_2)\leq F_n(\Theta_{01},\Theta_{02}).
\eqa
 That is
\bqa\label{convex}
 &&\sum_{k=1}^{2}\left\{tr(\widetilde{\Theta}_k-\Theta_{0k})\widehat{\Sigma}_k-\left(\log \det(\widetilde{\Theta}_k)-\log\det(\Theta_{0k})\right)
+\lambda||\widetilde{\Theta}^{-}_k||_{1}\right\}
+\rho||\widetilde{\Theta}_1^{-}-\widetilde{\Theta}_2^{-}||_1\ \non
&\leq& \lambda||\Theta^{-}_{01}||_{1}+\lambda||\Theta^{-}_{02}||_{1}
+\rho||\Theta_{01}^{-}-\Theta_{02}^{-}||_1.
\eqa
 Let $\Delta_k=\widetilde{\Theta}_k-\Theta_{0k}$, and
\begin{equation*}
f(\Delta_k):=tr(\Delta_k\Sigma_{0k})-\Big[\log\det(\Delta_k+\Theta_{0k})-\log\det(\Theta_{0k})\Big],
\end{equation*}
 subtracting $tr(\Delta_1(\widehat{\Sigma}_1-\Sigma_{01}))
+tr(\Delta_2(\widehat{\Sigma}_2-\Sigma_{02}))$ on the both sides of the inequality (\ref{convex}), we get
\begin{equation}\label{equ12}
\begin{aligned}
&f(\Delta_1)+f(\Delta_2)+\lambda||\widetilde{\Theta}_1^{-}||_1+\lambda||\widetilde{\Theta}_2^{-}||_1
+\rho||\widetilde{\Theta}_1^{-}-\widetilde{\Theta}_2^{-}||_1\\
\leq &-tr(\Delta_1(\widehat{\Sigma}_1-\Sigma_{01}))
-tr(\Delta_2(\widehat{\Sigma}_2-\Sigma_{02}))+
\lambda||\Theta_{01}^{-}||_1+\lambda||\Theta_{02}^{-}||_1
+\rho||\Theta_{01}^{-}-\Theta_{02}^{-}||_1.
\end{aligned}
\end{equation}

 For $tr(\Delta_k(\widehat{\Sigma}_k-\Sigma_{0k}))$ term, we have
\begin{equation}
\begin{aligned}
|tr(\Delta_k(\widehat{\Sigma}_k-\Sigma_{0k}))|&=|G(\Delta_k\circ (\widehat{\Sigma}_k-\Sigma_{0k}))|\\
&\leq|G(\Delta_k^{-}\circ (\widehat{\Sigma}_k^{-}-\Sigma_{0k}^{-}))|+|G(\Delta_k^{+}\circ (\widehat{\Sigma}_k^{+}-\Sigma_{0k}^{+}))|,\\
\end{aligned}
\end{equation}
 where function $G(M)$ takes the summation of all the elements of the matrix $M$, and $\circ$ is Hadamard product.
 According to Cauchy-Schwarz inequality, on the sets $\{\max_k\{||\widehat{\Sigma}_k-\Sigma_{0k}||_{\infty}\}\leq\lambda_0\}$,
\begin{equation}
\begin{aligned}
|G(\Delta_k^{-}\circ (\widehat{\Sigma}_k^{-}-\Sigma_{0k}^{-}))|+|G(\Delta_k^{+}\circ (\widehat{\Sigma}_k^{+}-\Sigma_{0k}^{+}))|
&\leq ||\widehat{\Sigma}_{k}^{-}-
\Sigma_{0k}^{-}||_{\infty}||\Delta_k^{-}||_1+||\widehat{\Sigma}_{k}^{+}-
\Sigma_{0k}^{+}||_F||\Delta_k^{+}||_{F}\\
&\leq \lambda_0||\Delta_k^{-}||_1+||\widehat{\Sigma}_{k}^{+}-
\Sigma_{0k}^{+}||_F||\Delta_k^{+}||_{F}.
\end{aligned}
\end{equation}
 Hence,
\begin{equation}\label{trabou}
\begin{aligned}
-tr(\Delta_k(\widehat{\Sigma}_k-\Sigma_{0k}))& \leq
|tr(\Delta_k(\widehat{\Sigma}_k-\Sigma_{0k}))|\\
&\leq \lambda_0||\Delta_k^{-}||_1+||\widehat{\Sigma}_{k}^{+}-
\Sigma_{0k}^{+}||_F||\Delta_k^{+}||_{F}.
\end{aligned}
\end{equation}

 Next, for $L_k\geq 1$ satisfying condition
\bqa
1/L_k\leq \lambda_{min}(\Theta_{0k})\leq \lambda_{max}(\Theta_{0k})\leq L_k,
\eqa
 we choose $L>1$ satisfying $1/L\leq 1/L_k$ and $L_k\leq L$, $k=1,2.$
 Based on the definitions of $\Delta_k$ and $\widetilde{\Theta}_k$, we get
\begin{equation}\label{deltaine}
\begin{aligned}
||\Delta_k||_{F}=\alpha_k||\widehat{\Theta}_k-\Theta_{0k}||_{F}
=\frac{||\widehat{\Theta}_k-\Theta_{0k}||_{F}}{M+||\widehat{\Theta}_k-\Theta_{0k}||_{F}}M,
\end{aligned}
\end{equation}
 for arbitrary $M$ in $(0,1/2L]$.
 Thus, $||\Delta_k||_{F}$ is bounded by $M$, i.e., $||\Delta_k||_{F}\leq M$.
 For $f(\Delta_k)$ term, based on Lemma \ref{lembasic}, we have
\begin{equation}\label{lembou}
f(\Delta_k) \geq c||\widetilde{\Theta}_k-\Theta_{0k}||_{F}^{2},
\end{equation}
 where $c=\frac{1}{2(L+1/(2L))^2}$. In particular, we choose $c=1/(8L^2)$, and the inequality (\ref{lembou}) still holds.

 Using bounds (\ref{trabou}) and (\ref{lembou}), the inequality (\ref{equ12}) turns to be
\begin{equation}\label{equatran}
\begin{aligned}
&c||\widetilde{\Theta}_1-\Theta_{01}||_{F}^{2}+
c||\widetilde{\Theta}_2-\Theta_{02}||_{F}^{2}+
\lambda||\widetilde{\Theta}_1^{-}||_1+
\lambda||\widetilde{\Theta}_2^{-}||_1+
\rho||\widetilde{\Theta}_1^{-}-\widetilde{\Theta}_2^{-}||_1\\
\leq &\lambda_0||\Delta_1^{-}||_1+\lambda_0||\Delta_2^{-}||_1
+||\widehat{\Sigma}_{1}^{+}-\Sigma_{01}^{+}||_F||\Delta_1^{+}||_{F}
+||\widehat{\Sigma}_{2}^{+}-\Sigma_{02}^{+}||_F||\Delta_2^{+}||_{F}\\
&+\lambda||\Theta_{01}^{-}||_1+\lambda||\Theta_{02}^{-}||_1
+\rho||\Theta_{01}^{-}-\Theta_{02}^{-}||_1.
\end{aligned}
\end{equation}
 We move some terms of the inequality (\ref{equatran}) and combine them to get the following inequality
\begin{equation}\label{equa3needpro}
\begin{split}
&c||\widetilde{\Theta}_1-\Theta_{01}||_{F}^{2}+c||\widetilde{\Theta}_2-\Theta_{02}||_{F}^{2} \\
&+\lambda\Big\{||\widetilde{\Theta}_1^{-}||_1-||\Theta_{01}^{-}||_1+||\widetilde{\Theta}_2^{-}||_1-||\Theta_{02}^{-}||_1\Big\}
\\
\leq &\lambda_0\Big\{||\widetilde{\Theta}_1^{-}-\Theta_{01}^{-}||_1+||\widetilde{\Theta}_2^{-}-\Theta_{02}^{-}||_1\Big\}
+\rho\Big\{||\Theta_{01}^{-}-\Theta_{02}^{-}||_1-||\widetilde{\Theta}_1^{-}-\widetilde{\Theta}_2^{-}||_1\Big\}\\
&+||\widehat{\Sigma}_{1}^{+}-\Sigma_{01}^{+}||_F||\widetilde{\Theta}_1^{+}-\Theta_{01}^{+}||_{F}
+||\widehat{\Sigma}_{2}^{+}-\Sigma_{02}^{+}||_F||\widetilde{\Theta}_2^{+}-\Theta_{02}^{+}||_{F}.
\end{split}
\end{equation}

 Next we need to prove three inequations:
\bqa\label{equa1}
&||\widetilde{\Theta}_k^{-}||_1-||\Theta_{0k}^{-}||_1 \geq ||\Delta_{kS_k^{c}}^{-}||_1-||\Delta_{kS_k}^{-}||_1,\\ \label{equa2}
&||\widetilde{\Theta}_k^{-}-\Theta_{0k}^{-}||_1\leq
||\Delta_{kS_k^{c}}^{-}||_1+||\Delta_{kS_k}^{-}||_1,\\ \label{equa3}
&||\Theta_{01}^{-}-\Theta_{02}^{-}||_1-
||\widetilde{\Theta}_1^{-}-\widetilde{\Theta}_2^{-}||_1\leq
||\widetilde{\Theta}_1^{-}-\Theta_{01}^{-}||_1+||\widetilde{\Theta}_2^{-}
-\Theta_{02}^{-}||_1.
\eqa
 Because
\begin{equation}
\begin{split}
||\widetilde{\Theta}_k^{-}||_1=&||{\Theta}_{0k}^{-}+\Delta_k^{-}||_1\\
=&||{\Theta}_{0kS_{k}}^{-}+\Delta_{kS_k}^{-}||_1+||\Delta_{kS_k^c}^{-}||_1,
\end{split}
\end{equation}
 and
\bqa
||{\Theta}_{0k}^{-}||_1=||{\Theta}_{0kS_{k}}^{-}||_1
\eqa
 hold. Thus,
\begin{equation}
\begin{split}
||\widetilde{\Theta}_k^{-}||_1-||{\Theta}_{0k}^{-}||_1
=&||{\Theta}_{0kS_{k}}^{-}
+\Delta_{kS_k}^{-}||_1+||\Delta_{kS_k^c}^{-}||_1-||{\Theta}_{0kS_{k}}^{-}||_1\\
\geq&||\Delta_{kS_k^c}^{-}||_1-\Big|||{\Theta}_{0kS_{k}}^{-}
+\Delta_{kS_k}^{-}||_1-||{\Theta}_{0kS_{k}}^{-}||_1\Big|\\
\geq&||\Delta_{kS_k^c}^{-}||_1-||\Delta_{kS_k}^{-}||_1,
\end{split}
\end{equation}
 which proves inequality (\ref{equa1}).
 By the triangle inequality, we naturally obtain
\begin{equation}
\begin{split}
||\widetilde{\Theta}_k^{-}-\Theta_{0k}^{-}||_1=&||\Delta_k^{-}||_1\\
=&||\Delta_{kS_k^{c}}^{-}+\Delta_{kS_k}^{-}||_1\\
\leq& ||\Delta_{kS_k^{c}}^{-}||_1+||\Delta_{kS_k}^{-}||_1.
\end{split}
\end{equation}
 Thus, the inequation (\ref{equa2}) holds.
 For inequation (\ref{equa3}), we have
\begin{equation}
\begin{split}
||\Theta_{01}^{-}-\Theta_{02}^{-}||_1-||\widetilde{\Theta}_1^{-}-\widetilde{\Theta}_2^{-}||_1
=&||\Theta_{01}^{-}-\widetilde{\Theta}_1^{-}+\widetilde{\Theta}_1^{-}-\widetilde{\Theta}_2^{-}
+\widetilde{\Theta}_2^{-}
-\Theta_{02}^{-}||_1-||\widetilde{\Theta}_1^{-}-\widetilde{\Theta}_2^{-}||_1\\
\leq& ||\widetilde{\Theta}_1^{-}-\Theta_{01}^{-}||_1+||\widetilde{\Theta}_2^{-}
-\Theta_{02}^{-}||_1.
\end{split}
\end{equation}
 Thus, the inequality (\ref{equa3needpro}) yields
\begin{equation}
\begin{split}
&c||\widetilde{\Theta}_1-\Theta_{01}||_{F}^{2}+c||\widetilde{\Theta}_2-\Theta_{02}||_{F}^{2} \\
&+\lambda\Big\{||\Delta_{1S_1^{c}}^{-}||_1-||\Delta_{1S_1}^{-}||_1+||\Delta_{2S_2^{c}}^{-}||_1-||\Delta_{2S_2}^{-}||_1\Big\}\\
\leq &(\rho+\lambda_0)\Big\{||\Delta_{1S_1^{c}}^{-}||_1+||\Delta_{1S_1}^{-}||_1+||\Delta_{2S_2^{c}}^{-}||_1+||\Delta_{2S_2}^{-}||_1\Big\}\\
&+||\widehat{\Sigma}_{1}^{+}-\Sigma_{01}^{+}||_F||\widetilde{\Theta}_1^{+}-\Theta_{01}^{+}||_{F}
+||\widehat{\Sigma}_{2}^{+}-\Sigma_{02}^{+}||_F||\widetilde{\Theta}_2^{+}-\Theta_{02}^{+}||_{F}.
\end{split}
\end{equation}

 By taking $2(\rho+\lambda_0)<\lambda$, we conclude that
\begin{equation}
\begin{split}
&2c\Big\{||\widetilde{\Theta}_1-\Theta_{01}||_{F}^{2}+||\widetilde{\Theta}_2-\Theta_{02}||_{F}^{2}\Big\}
+\lambda\Big\{||\Delta_{1S_1^{c}}^{-}||_1+||\Delta_{2S_2^{c}}^{-}||_1\Big\}\\
\leq &3\lambda\Big\{||\Delta_{1S_1}^{-}||_1+||\Delta_{2S_2}^{-}||_1\Big\}\\
&+2\Big\{||\widehat{\Sigma}_{1}^{+}-\Sigma_{01}^{+}||_F||\widetilde{\Theta}_1^{+}-\Theta_{01}^{+}||_{F}
+||\widehat{\Sigma}_{2}^{+}-\Sigma_{02}^{+}||_F||\widetilde{\Theta}_2^{+}-\Theta_{02}^{+}||_{F}\Big\}.
\end{split}
\end{equation}
 By the definition of $\Delta_k$, we have
\begin{equation}\label{equal_35}
\begin{split}
||\Delta_k^{-}||_1=&||\Delta_{kS_{k}}^{-}+\Delta_{kS_{k}^c}^{-}||_1\\
\leq& ||\Delta_{kS_{k}}^{-}||_1+||\Delta_{kS_{k}^c}^{-}||_1.
\end{split}
\end{equation}
 So we deduce
\begin{equation}
\begin{split}
&2c\Big\{||\widetilde{\Theta}_1-\Theta_{01}||_{F}^{2}+||\widetilde{\Theta}_2-\Theta_{02}||_{F}^{2}\Big\}
+\lambda\Big\{||\Delta_{1}^{-}||_1+||\Delta_{2}^{-}||_1\Big\}\\
\leq &4\lambda\Big\{||\Delta_{1S_1}^{-}||_1+||\Delta_{2S_2}^{-}||_1\Big\}\\
&+2\Big\{||\widehat{\Sigma}_{1}^{+}-\Sigma_{01}^{+}||_F||\widetilde{\Theta}_1^{+}-\Theta_{01}^{+}||_{F}
+||\widehat{\Sigma}_{2}^{+}-\Sigma_{02}^{+}||_F||\widetilde{\Theta}_2^{+}-\Theta_{02}^{+}||_{F}\Big\}
\end{split}
\end{equation}
 holds.
 Since the inequality of arithmetic and geometric means,
% $(|x_1|+|x_2|+...+|x_n|)^2\leq n(|x_1|^2+|x_2|^2+...+|x_n|^2)$,$xy\leq \frac{1}{2}(x^2+y^2)$
 the inequality $||\Delta_{kS_k}^{-}||_1\leq \sqrt{s_k}||\Delta_{kS_k}^{-}||_{F}$ holds. Thus
\begin{equation}\label{lasine}
\begin{split}
&2c\Big\{||\widetilde{\Theta}_1-\Theta_{01}||_{F}^{2}+||\widetilde{\Theta}_2-\Theta_{02}||_{F}^{2}\Big\}
+\lambda\Big\{||\Delta_{1}^{-}||_1+||\Delta_{2}^{-}||_1\Big\}\\
\leq &4\lambda\Big\{\sqrt{s_1}||\Delta_{1S_1}^{-}||_{F}+\sqrt{s_2}||\Delta_{2S_2}^{-}||_{F}\Big\}\\
&+2\Big\{||\widehat{\Sigma}_{1}^{+}-\Sigma_{01}^{+}||_F||\widetilde{\Theta}_1^{+}-\Theta_{01}^{+}||_{F}
+||\widehat{\Sigma}_{2}^{+}-\Sigma_{02}^{+}||_F||\widetilde{\Theta}_2^{+}-\Theta_{02}^{+}||_{F}\Big\}.\\
\end{split}
\end{equation}
 Using $xy\leq (x^2+y^2)/2$, the inequality (\ref{lasine}) infer that
\begin{equation}
\begin{split}
&2c\Big\{||\widetilde{\Theta}_1-\Theta_{01}||_{F}^{2}+||\widetilde{\Theta}_2-\Theta_{02}||_{F}^{2}\Big\}
+\lambda\Big\{||\Delta_{1}^{-}||_1+||\Delta_{2}^{-}||_1\Big\}\\
\leq
 &\frac{1}{2}\Big(c||\Delta_{1S_1}^{-}||_{F}^{2}+\frac{16\lambda^2s_1}{c}
 +c||\Delta_{2S_2}^{-}||_{F}^{2}+\frac{16\lambda^2s_2}{c}\Big)\\
&+\frac{1}{2}\Big(c||\widetilde{\Theta}_1^{+}-\Theta_{01}^{+}||_{F}^{2}+
\frac{4||\widehat{\Sigma}_{1}^{+}-\Sigma_{01}^{+}||_{F}^{2}}{c}
+c||\widetilde{\Theta}_2^{+}-\Theta_{02}^{+}||_{F}^{2}+
\frac{4||\widehat{\Sigma}_{2}^{+}-\Sigma_{02}^{+}||_{F}^{2}}{c}\Big).\\
\end{split}
\end{equation}
 Because
\begin{equation}\label{equal_39}
\begin{split}
c||\widetilde{\Theta}_k^{+}-\Theta_{0k}^{+}||_{F}^{2}+c||\Delta_{kS_k}^{-}||_{F}^{2}
\leq
&\Big\{c||\widetilde{\Theta}_k^{+}-\Theta_{0k}^{+}||_{F}^{2}+c||\Delta_{k}^{-}||_{F}^{2}\Big\}\\
&+\Big\{c||\Delta_{kS_k}^{-}||_{F}^{2}+c||\Delta_{kS_k^{c}}^{-}||_{F}^{2}+c||\Delta_{k}^{+}||_{F}^{2}\Big\}\\
=&2c||\Delta_{k}||_{F}^{2},
\end{split}
\end{equation}
 we obtain
\begin{equation}
\begin{split}
&2c\Big\{||\widetilde{\Theta}_1-\Theta_{01}||_{F}^{2}+||\widetilde{\Theta}_2-\Theta_{02}||_{F}^{2}\Big\}
+\lambda\Big\{||\Delta_{1}^{-}||_1+||\Delta_{2}^{-}||_1\Big\}\\
\leq &c\Big\{||\Delta_{1}||_{F}^{2}+||\Delta_{2}||_{F}^{2}\Big\}+\frac{8\lambda^2(s_1+s_2)}{c}
+\frac{2||\widehat{\Sigma}_{1}^{+}-\Sigma_{01}^{+}||_{F}^{2}}{c}
+\frac{2||\widehat{\Sigma}_{2}^{+}-\Sigma_{02}^{+}||_{F}^{2}}{c}.
\end{split}
\end{equation}
 Thus,
\begin{equation}\label{equbf40}
\begin{split}
&c\Big\{||\Delta_{1}||_{F}^{2}+||\Delta_{2}||_{F}^{2}\Big\}
+\lambda\Big\{||\Delta_{1}^{-}||_1+||\Delta_{2}^{-}||_1\Big\}\\
\leq &\frac{8\lambda^2(s_1+s_2)}{c}
+\frac{2||\widehat{\Sigma}_{1}^{+}-\Sigma_{01}^{+}||_{F}^{2}}{c}
+\frac{2||\widehat{\Sigma}_{2}^{+}-\Sigma_{02}^{+}||_{F}^{2}}{c}.
\end{split}
\end{equation}
 Based on the inequality $||\widehat{\Sigma}_{k}^{+}-\Sigma_{0k}^{+}||_{F}\leq \sqrt{p}||\widehat{\Sigma}_{k}^{+}-\Sigma_{0k}^{+}||_{\infty}$, we have
\begin{equation}\label{oracleinqua}
\begin{split}
c\Big\{||\Delta_{1}||_{F}^{2}+||\Delta_{2}||_{F}^{2}\Big\}
+\lambda\Big\{||\Delta_{1}^{-}||_1+||\Delta_{2}^{-}||_1\Big\}
\leq  \frac{8\lambda^2(s_1+s_2)}{c}+\frac{4p\lambda_0^2}{c}.
\end{split}
\end{equation}

 Next, we prove that substituting $\widehat{\Theta}_k$ for $\widetilde{\Theta}_k$, the conclusion still holds.
 According to the condition,
\begin{equation}
\begin{split}
||\Delta_{1}||_{F}^{2}+||\Delta_{2}||_{F}^{2}
\leq &\frac{\lambda_0}{2cL} \leq \frac{\lambda}{4cL}\leq\frac{1}{32L^2}.
\end{split}
\end{equation}
 Taking $M=1/(2\sqrt{2}L)<1/2L$, we have
\begin{equation}
\begin{split}
||\Delta_{1}||_{F}^{2}+||\Delta_{2}||_{F}^{2}
\leq &M^2/4.
\end{split}
\end{equation}
 Thus, $||\Delta_{k}||_{F}$ is bounded by $M/2$.
 In addition,
\begin{equation}
\begin{split}
||\widehat{\Theta}_k-\Theta_{0k}||_F=\frac{M||\Delta_k||_F}{M-||\Delta_k||_F},
\end{split}
\end{equation}
 which means $||\widehat{\Theta}_k-\Theta_{0k}||_F$ is monotone increasing function of $||\Delta_k||_F$ on set $(0,M)$. We obtain that $||\widehat{\Theta}_k-\Theta_{0k}||_F\leq M$. Therefore, we can substitute $\widehat{\Theta}_k$ for $\widetilde{\Theta}_k$, and that leads to the inequality (\ref{oracleinqua}) holds for $\widehat{\Theta}_k$.

 According to inequality (\ref{oracleinqua}), we get
\begin{equation}
\begin{split}
||\widehat{\Theta}_k-\Theta_{0k}||_{F}^2\leq& \frac{8\lambda^2(s_1+s_2)}{c^2}+\frac{4p\lambda_0^2}{c^2}\\
\leq &\frac{\lambda^2(8s_1+8s_2+p)}{c^2},
\end{split}
\end{equation}
 and
\begin{equation}
\begin{split}
||\widehat{\Theta}_k^{-}-\Theta_{0k}^{-}||_{1}\leq& \frac{8\lambda(s_1+s_2)}{c}+\frac{4p\lambda_0^2}{\lambda c}\\
\leq &\frac{\lambda(8s_1+8s_2+p)}{c}.
\end{split}
\end{equation}
  Thus, we conclude the upper bound of $\sum_{k=1}^2|||\widehat{\Theta}_k-\Theta_{0k}|||_1$,
\begin{equation}
\begin{split}
\sum_{k=1}^2|||\widehat{\Theta}_k-\Theta_{0k}|||_1
\leq &
\sum_{k=1}^2\Big(||\widehat{\Theta}_k^{+}-\Theta_{0k}^{+}||_{\infty}
+||\widehat{\Theta}_k^{-}-\Theta_{0k}^{-}||_{1}\Big)\\
\leq &\sum_{k=1}^2\Big(||\widehat{\Theta}_k-\Theta_{0k}||_{F}
+||\widehat{\Theta}_k^{-}-\Theta_{0k}^{-}||_{1}\Big)\\
\leq &\frac{2\lambda\sqrt{8s_1+8s_2+p}}{c}+\frac{2\lambda(8s_1+8s_2+p)}{c}\\
\leq &\frac{4\lambda(8s_1+8s_2+p)}{c}.
\end{split}
\end{equation}

\end{proof}

\subsection{Proof of Theorem \ref{weiandr}}

\begin{proof}
 The minimizer $(\widehat{\Theta}_R^{[1]},\widehat{\Theta}_R^{[2]})$ satisfying inequality (\ref{equbf40}), that is
\begin{equation}
\begin{split}
&c\Big\{||\widehat{\Theta}_R^{[1]}-\Theta_{R0}^{[1]}||_{F}^{2}+
||\widehat{\Theta}_R^{[2]}-\Theta_{R0}^{[2]}||_{F}^{2}\Big\}
+\lambda\Big\{||(\widehat{\Theta}_R^{[1]}-\Theta_{R0}^{[1]})^{-}||_1
+||(\widehat{\Theta}_R^{[2]}-\Theta_{R0}^{[2]})^{-}||_1\Big\}\\
\leq &\frac{8\lambda^2(s_1+s_2)}{c}
+\frac{2||(\widehat{R}^{[1]}-R_0^{[1]})^{+}||_{F}^{2}}{c}
+\frac{2||(\widehat{R}^{[2]}-R_0^{[2]})^{+}||_{F}^{2}}{c}.
\end{split}
\end{equation}
 The diagonal elements of $\widehat{R}^{[k]}$ and $R_{0}^{[k]}$ are all $1$. Thus
\begin{equation}
\begin{split}
&c\Big\{||\widehat{\Theta}_R^{[1]}-\Theta_{R0}^{[1]}||_{F}^{2}+
||\widehat{\Theta}_R^{[2]}-\Theta_{R0}^{[2]}||_{F}^{2}\Big\}
+\lambda\Big\{||(\widehat{\Theta}_R^{[1]}-\Theta_{R0}^{[1]})^{-}||_1
+||(\widehat{\Theta}_R^{[2]}-\Theta_{R0}^{[2]})^{-}||_1\Big\}\\
\leq &\frac{8\lambda^2(s_1+s_2)}{c}.
\end{split}
\end{equation}

 Moreover, for the conclusion of the $l_1$-operator norm, we get
\begin{equation}\label{bou}
\begin{split}
|||\widehat{\Theta}_R^{[1]}-\Theta_{R0}^{[1]}|||_1
+|||\widehat{\Theta}_R^{[2]}-\Theta_{R0}^{[2]}|||_1
\leq &
\sum_{k=1}^2\Big(||(\widehat{\Theta}_R^{[k]}-\Theta_{R0}^{[k]})^{+}||_{\infty}
+||(\widehat{\Theta}_R^{[k]}-\Theta_{R0}^{[k]})^{-}||_{1}\Big)\\
\leq &\sum_{k=1}^2\Big(||\widehat{\Theta}_R^{[k]}-\Theta_{R0}^{[k]}||_{F}
+||(\widehat{\Theta}_R^{[k]}-\Theta_{R0}^{[k]})^{-}||_{1}\Big)\\
\leq &\frac{32\lambda(s_1+s_2)}{c}.
\end{split}
\end{equation}

 For the minimizer $(\widehat{\Theta}_w^{[1]},\widehat{\Theta}_w^{[2]})$, following inequality holds
\begin{equation}\label{weiineq}
\begin{split}
|||\widehat{\Theta}_R^{[k]}-\Theta_{R0}^{[k]}|||_1
=&
|||\widehat{W}^{[k]}\widehat{\Theta}_w^{[k]}\widehat{W}^{[k]}
-W_0^{[k]}\Theta_{w0}^{[k]}W_0^{[k]}|||_1\\
\leq &||\widehat{W}^{[k]}||_{\infty}^2||| \widehat{\Theta}_w^{[k]}- \Theta_{w0}^{[k]} |||_1
+||  \widehat{W}^{[k]}- W_0^{[k]}  ||_{\infty}||| \Theta_{w0}^{[k]}  |||_1
|| \widehat{W}^{[k]} ||_{\infty}\\
&+|| W_0^{[k]}  ||_{\infty} |||  \Theta_{w0}^{[k]}  |||_1||  \widehat{W}^{[k]}- W_0^{[k]} ||_{\infty}.
\end{split}
\end{equation}
 To draw the conclusion, we have the following facts:
 \begin{itemize}
   \item The Sub-Gaussian vector with covariance $\Sigma_0^{[k]}$ implies that $\sqrt{n/\log{p}}||  (\widehat{\Sigma}^{[k]}-\Sigma_0^{[k]}) ||_{\infty}$ is bounded in probability.
   \item The eigenvalues of $\Theta_{w0}^{[k]}$ are bounded by a constant.
 \end{itemize}
  Thus, $|||\widehat{\Theta}_R^{[k]}-\Theta_{R0}^{[k]}|||_1$ and $||| \widehat{\Theta}_w^{[k]}- \Theta_{w0}^{[k]} |||_1$ share the same boundary.

\end{proof}

\subsection{Proof of Theorem \ref{FGL}}

\begin{proof}

 Similarly, $\widehat{\Theta}_k$ are the minimum value of the fused graphical Lasso for $k=1,2,\cdots,K$. Let $\widetilde{\Theta}_k=\alpha_k\widehat{\Theta}_k+(1-\alpha_k)\Theta_{0k}$, and $\alpha_k=\frac{M}{M+||\widehat{\Theta}_k-\Theta_{0k}||_{F}}$.
 Denotes
 $$F_n(\Theta_1,\cdots,\Theta_K)=\sum_{k=1}^{K}\left\{tr(\Theta_k\widehat{\Sigma}_k)-\log\det(\Theta_k)\right\}
 +\lambda\sum_{k=1}^{K}||\Theta_k^{-}||_1+\rho\sum_{k<k'}||\Theta_k^{-}-\Theta_{k'}^{-}||_1,$$ we obtain
\bqa
F_n(\widetilde\Theta_1,\widetilde\Theta_2,\cdots,\widetilde\Theta_K)\leq F_n(\Theta_{01},\Theta_{02},\cdots,\Theta_{0K}).
\eqa
 Thus,
\bqa\label{convexFGL}
&&\sum_{k=1}^{K}\left\{tr(\widetilde{\Theta}_k-\Theta_{0k})\widehat{\Sigma}_k-\left(\log \det(\widetilde{\Theta}_k)-\log\det(\Theta_{0k})\right)
+\lambda||\widetilde{\Theta}^{-}_k||_{1}\right\}
+\rho\sum_{k<k'}||\widetilde{\Theta}_k^{-}-\widetilde{\Theta}_{k'}^{-}||_1\ \non
&\leq& \lambda\sum_{k=1}^{K}||\Theta^{-}_{0k}||_{1}
+\rho\sum_{k<k'}||\Theta_{0k}^{-}-\Theta_{0k'}^{-}||_1.
\eqa
 Using the notations that $\Delta_k=\widetilde{\Theta}_k-\Theta_{0k}$ and $$f(\Delta_k):=tr(\Delta_k\Sigma_{0k})-\Big[\log\det(\Delta_k+\Theta_{0k})-\log\det(\Theta_{0k})\Big]$$ we yield the following expression
\begin{equation}\label{equ12FGL}
\begin{aligned}
&\sum_{k=1}^{K}f(\Delta_k)+\lambda\sum_{k=1}^{K}||\widetilde{\Theta}_k^{-}||_1
+\rho\sum_{k<k'}||\widetilde{\Theta}_k^{-}-\widetilde{\Theta}_{k'}^{-}||_1\\
\leq &-\sum_{k=1}^{K}\left(tr(\Delta_k(\widehat{\Sigma}_k-\Sigma_{0k}))\right)
-tr(\Delta_2(\widehat{\Sigma}_2-\Sigma_{02}))+
\lambda\sum_{k=1}^{K}||\Theta_{0k}^{-}||_1
+\rho\sum_{k<k'}||\Theta_{0k}^{-}-\Theta_{0k'}^{-}||_1.
\end{aligned}
\end{equation}

 For $L_k\geq 1,k=1,2,\cdots,K$, the minimum and maximum eigenvalues of $\Theta_{0k}$ hold that
\bqa
1/L_k\leq \lambda_{min}(\Theta_{0k})\leq \lambda_{max}(\Theta_{0k})\leq L_k.
\eqa
 For multiple case, we select a constant $L$ satisfying $1/L\leq 1/L_k$ and $L_k\leq L$.
 By similar analysis, for $M$ in $(0,1/2L]$, the inequality (\ref{deltaine}) and the inequality (\ref{lembou}) still hold.

 For $K$ groups data, based on the inequalities (\ref{trabou}) and (\ref{lembou}).
 Then, the inequality (\ref{equ12FGL}) turns to be
\begin{equation}
\begin{aligned}
&c\sum_{k=1}^{K}||\widetilde{\Theta}_k-\Theta_{0k}||_{F}^{2}+
\lambda\sum_{k=1}^{K}||\widetilde{\Theta}_k^{-}||_1+
\rho\sum_{k<k'}||\widetilde{\Theta}_k^{-}-\widetilde{\Theta}_{k'}^{-}||_1\\
\leq &\sum_{k=1}^{K}\left\{\lambda_0||\Delta_k^{-}||_1
+||\widehat{\Sigma}_{k}^{+}-\Sigma_{0k}^{+}||_F||\Delta_k^{+}||_{F}\right\}+\lambda\sum_{k=1}^{K}||\Theta_{0k}^{-}||_1
+\rho\sum_{k<k'}||\Theta_{0k}^{-}-\Theta_{0k'}^{-}||_1.
\end{aligned}
\end{equation}
 Thus,
\begin{equation}\label{equa3needpro FGL}
\begin{split}
&c\sum_{k=1}^{K}||\widetilde{\Theta}_k-\Theta_{0k}||_{F}^{2}
+\lambda\sum_{k=1}^{K}\left\{||\widetilde{\Theta}_k^{-}||_1-||\Theta_{0k}^{-}||_1\right\}\\
\leq &\lambda_0\sum_{k=1}^{K}||\widetilde{\Theta}_k^{-}-\Theta_{0k}^{-}||_1
+\rho\sum_{k<k'}\Big\{||\Theta_{0k}^{-}-\Theta_{0k'}^{-}||_1-||\widetilde{\Theta}_k^{-}-\widetilde{\Theta}_{k'}^{-}||_1\Big\}\\
&+\sum_{k=1}^{K}\left\{||\widehat{\Sigma}_{k}^{+}-\Sigma_{0k}^{+}||_F||\widetilde{\Theta}_k^{+}-\Theta_{0k}^{+}||_{F}\right\}.
\end{split}
\end{equation}

 When $k=1,2,\cdots,K$, the inequations (\ref{equa1}) and (\ref{equa2}) still hold. Similarly, we have the following inequality
\begin{equation}\label{equa3FGL}
\begin{split}
||\Theta_{0k}^{-}-\Theta_{0k'}^{-}||_1-||\widetilde{\Theta}_k^{-}-\widetilde{\Theta}_{k'}^{-}||_1
=&||\Theta_{0k}^{-}-\widetilde{\Theta}_k^{-}+\widetilde{\Theta}_k^{-}-\widetilde{\Theta}_{k'}^{-}
+\widetilde{\Theta}_{k'}^{-}
-\Theta_{0k'}^{-}||_1-||\widetilde{\Theta}_k^{-}-\widetilde{\Theta}_{k'}^{-}||_1\\
\leq& ||\widetilde{\Theta}_k^{-}-\Theta_{0k}^{-}||_1+||\widetilde{\Theta}_{k'}^{-}
-\Theta_{0k'}^{-}||_1.
\end{split}
\end{equation}
Thus, by the equations (\ref{equa1}), (\ref{equa2}) and (\ref{equa3FGL}) the inequality (\ref{equa3needpro FGL}) yields
\begin{equation}
\begin{split}
&c\sum_{k=1}^{K}||\widetilde{\Theta}_k-\Theta_{0k}||_{F}^{2}
+\lambda\sum_{k=1}^{K}\Big\{||\Delta_{kS_k^{c}}^{-}||_1-||\Delta_{kS_k}^{-}||_1\Big\}\\
\leq &\lambda_0\sum_{k=1}^{K}\Big\{||\Delta_{kS_k^{c}}^{-}||_1+||\Delta_{kS_k}^{-}||_1\Big\}
+\rho\sum_{k<k'}\Big\{||\Delta_{kS_k^{c}}^{-}||_1+||\Delta_{kS_k}^{-}||_1+||\Delta_{k'S_{k'}^{c}}^{-}||_1+||\Delta_{k'S_{k'}}^{-}||_1\Big\}\\
&+\sum_{k=1}^{K}\left\{||\widehat{\Sigma}_{k}^{+}-\Sigma_{0k}^{+}||_F||\widetilde{\Theta}_k^{+}-\Theta_{0k}^{+}||_{F}\right\}\\
\leq &\left(\frac{K(K-1)}{2}\rho+\lambda_0\right)\sum_{k=1}^{K}\Big\{||\Delta_{kS_k^{c}}^{-}||_1+||\Delta_{kS_k}^{-}||_1\Big\}
+\sum_{k=1}^{K}\left\{||\widehat{\Sigma}_{k}^{+}-\Sigma_{0k}^{+}||_F||\widetilde{\Theta}_k^{+}-\Theta_{0k}^{+}||_{F}\right\}.
\end{split}
\end{equation}

 Since $K$ is a fixed constant, and $2(\frac{K(K-1)}{2}\rho+\lambda_0)<\lambda$, we can obtain
\begin{equation}
\begin{split}
&2c\sum_{k=1}^{K}||\widetilde{\Theta}_k-\Theta_{0k}||_{F}^{2}
+\lambda\sum_{k=1}^{K}||\Delta_{kS_k^{c}}^{-}||_1\\
\leq &3\lambda\sum_{k=1}^{K}||\Delta_{kS_k}^{-}||_1
+2\sum_{k=1}^{K}\Big\{||\widehat{\Sigma}_{k}^{+}-\Sigma_{0k}^{+}||_F||\widetilde{\Theta}_k^{+}-\Theta_{0k}^{+}||_{F}\Big\}.
\end{split}
\end{equation}

On the basis of the inequality (\ref{equal_35}), we deduce
\begin{equation}\label{Begin}
\begin{split}
&2c\sum_{k=1}^{K}||\widetilde{\Theta}_k-\Theta_{0k}||_{F}^{2}
+\lambda\sum_{k=1}^{K}||\Delta_{k}^{-}||_1\\
\leq &4\lambda\sum_{k=1}^{K}||\Delta_{kS_k}^{-}||_1
+2\sum_{k=1}^{K}\Big\{||\widehat{\Sigma}_{k}^{+}-\Sigma_{0k}^{+}||_F||\widetilde{\Theta}_k^{+}-\Theta_{0k}^{+}||_{F}\Big\}
\end{split}
\end{equation}
 holds.
 In addition, one can get the inequality $||\Delta_{kS_k}^{-}||_1\leq \sqrt{s_k}||\Delta_{kS_k}^{-}||_{F}$. Thus
\begin{equation}\label{lasineFGL}
\begin{split}
&2c\sum_{k=1}^{K}||\widetilde{\Theta}_k-\Theta_{0k}||_{F}^{2}
+\lambda\sum_{k=1}^{K}||\Delta_{k}^{-}||_1\\
\leq &4\lambda\sum_{k=1}^{K}\left(\sqrt{s_k}||\Delta_{kS_k}^{-}||_{F}\right)
+2\sum_{k=1}^{K}\Big\{||\widehat{\Sigma}_{k}^{+}-\Sigma_{0k}^{+}||_F||\widetilde{\Theta}_k^{+}-\Theta_{0k}^{+}||_{F}\Big\}.
\end{split}
\end{equation}
 Based on $xy\leq (x^2+y^2)/2$ and the inequality (\ref{equal_39}), the inequality (\ref{lasineFGL}) infer that
\begin{equation}
\begin{split}
&2c\sum_{k=1}^{K}||\widetilde{\Theta}_k-\Theta_{0k}||_{F}^{2}
+\lambda\sum_{k=1}^{K}||\Delta_{k}^{-}||_1\\
\leq
&\frac{1}{2}\sum_{k=1}^{K}\Big(c||\Delta_{kS_k}^{-}||_{F}^{2}+\frac{16\lambda^2s_k}{c}\Big)
+\frac{1}{2}\sum_{k=1}^{K}\Big(c||\widetilde{\Theta}_k^{+}-\Theta_{0k}^{+}||_{F}^{2}+
\frac{4||\widehat{\Sigma}_{k}^{+}-\Sigma_{0k}^{+}||_{F}^{2}}{c}\Big)\\
\leq
&c\sum_{k=1}^{K}||\Delta_{k}||_{F}^{2}+\frac{8\lambda^2\sum_{k=1}^{K}s_k}{c}
+\frac{2\sum_{k=1}^{K}||\widehat{\Sigma}_{k}^{+}-\Sigma_{0k}^{+}||_{F}^{2}}{c}.
\end{split}
\end{equation}
Thus,
\begin{equation}\label{equbf40FGL}
c\sum_{k=1}^{K}||\Delta_{k}||_{F}^{2}+\lambda\sum_{k=1}^{K}||\Delta_{k}^{-}||_1
\leq
\frac{8\lambda^2\sum_{k=1}^{K}s_k}{c}
+\frac{2\sum_{k=1}^{K}||\widehat{\Sigma}_{k}^{+}-\Sigma_{0k}^{+}||_{F}^{2}}{c}.
\end{equation}
Using the relation between the Frobenius norm and the supremum norm, we have
\begin{equation}\label{oracleinqua FGL}
\begin{split}
c\sum_{k=1}^{K}||\Delta_{k}||_{F}^{2}+\lambda\sum_{k=1}^{K}||\Delta_{k}^{-}||_1
\leq  \frac{8\lambda^2\sum_{k=1}^{K}s_k}{c}+\frac{2Kp\lambda_0^2}{c}.
\end{split}
\end{equation}

 According to the inequality (\ref{oracleinqua FGL}), we get
\bqa
\sum_{k=1}^{K}||\Delta_{k}||_{F}^{2}
\leq \frac{\lambda_0}{2cL}.
\eqa
 According to $\lambda_0\leq \lambda/2$ and the condition $\lambda\leq c/8L$, we get
\begin{equation}
\sum_{k=1}^{K}||\Delta_{k}||_{F}^{2}
\leq\frac{1}{32L^2}.
\end{equation}
 Taking $M=1/(2\sqrt{2}L)<1/2L$, we have
\begin{equation}
\sum_{k=1}^{K}||\Delta_{k}||_{F}^{2}
\leq M^2/4.
\end{equation}
Thus, $||\Delta_{k}||_{F}$ is bounded by $M/2$.
Further, we can derive $||\widehat{\Theta}_k-\Theta_{0k}||_F\leq M$ which means that we can substitute $\widehat{\Theta}_k$ for $\widetilde{\Theta}_k$, and that leads to the inequality (\ref{oracleinqua FGL}) holds for $\widehat{\Theta}_k$, i.e.
\begin{equation}
\begin{split}
c\sum_{k=1}^{K}||\widehat{\Theta}_k-\Theta_{0k}||_{F}^2
+\lambda\sum_{k=1}^{K}||(\widehat{\Theta}_k-\Theta_{0k})^{-}||_1
\leq \frac{8\lambda^2\sum_{k=1}^{K}s_k}{c}+\frac{2Kp\lambda_0^2}{c},
\end{split}
\end{equation}
That implies
\begin{equation}\label{End}
\begin{split}
\sum_{k=1}^K|||\widehat{\Theta}_k-\Theta_{0k}|||_1
\leq &
\sum_{k=1}^K\Big(||\widehat{\Theta}_k^{+}-\Theta_{0k}^{+}||_{\infty}
+||\widehat{\Theta}_k^{-}-\Theta_{0k}^{-}||_{1}\Big)\\
\leq &\sum_{k=1}^K\Big(||\widehat{\Theta}_k-\Theta_{0k}||_{F}
+||\widehat{\Theta}_k^{-}-\Theta_{0k}^{-}||_{1}\Big)\\
\leq &K\left[\frac{\lambda\sqrt{8\sum_{k=1}^{K}s_k+\frac{Kp}{2}}}{c}
+\frac{\lambda\left(8\sum_{k=1}^{K}s_k+\frac{Kp}{2}\right)}{c}\right]\\
\leq &\frac{2K\lambda\left(8\sum_{k=1}^{K}s_k+\frac{Kp}{2}\right)}{c},
\end{split}
\end{equation}
 which completes the proof.

\end{proof}

\subsection{Proof of Theorem \ref{FGLwei}}

\begin{proof}

 We get from (\ref{equbf40FGL})
\begin{equation}
c\sum_{k=1}^{K}||\widehat{\Theta}_R^{[k]}-\Theta_{R0}^{[k]}||_{F}^{2}
+\lambda\sum_{k=1}^{K}||(\widehat{\Theta}_R^{[k]}-\Theta_{R0}^{[k]})^{-}||_1
\leq
\frac{8\lambda^2\sum_{k=1}^{K}s_k}{c}
+\frac{2\sum_{k=1}^{K}||(\widehat{\Theta}_R^{[k]}-\Theta_{R0}^{[k]})^{+}||_{F}^{2}}{c},
\end{equation}
 and similarly derive
\begin{equation}
c\sum_{k=1}^{K}||\widehat{\Theta}_R^{[k]}-\Theta_{R0}^{[k]}||_{F}^{2}
+\lambda\sum_{k=1}^{K}||(\widehat{\Theta}_R^{[k]}-\Theta_{R0}^{[k]})^{-}||_1
\leq
\frac{8\lambda^2\sum_{k=1}^{K}s_k}{c}.
\end{equation}

 Using
\begin{equation}
\sum_{k=1}^{K}|||\widehat{\Theta}_R^{[k]}-\Theta_{R0}^{[k]}|||_{1}\leq
\sum_{k=1}^{K}\left(||\widehat{\Theta}_R^{[k]}-\Theta_{R0}^{[k]}||_{F}
+||(\widehat{\Theta}_R^{[k]}-\Theta_{R0}^{[k]})^{-}||_1\right)
\end{equation}
 we have
\begin{equation}
\sum_{k=1}^{K}|||\widehat{\Theta}_R^{[k]}-\Theta_{R0}^{[k]}|||_{1}
\leq
\frac{16K\lambda\sum_{k=1}^{K}s_k}{c}.
\end{equation}

 At last, using the inequality (\ref{weiineq}), based on the analysis of the upper bound of $|| W_0^{[k]}  ||_{\infty}$ and $|| \widehat{W}^{[k]} ||_{\infty}$, and the convergence rate of $||  (\widehat{\Sigma}^{[k]}-\Sigma_0^{[k]}) ||_{\infty}$, we draw the conclusion that
\begin{equation}
\sum_{k=1}^{K}|||\widehat{\Theta}_w^{[k]}-\Theta_{0}^{[k]}|||_1
\frac{16K\lambda\sum_{k=1}^{K}s_k}{c}.
\end{equation}

\end{proof}

\subsection{Proof of Theorem \ref{CLT}}

\begin{proof}
 First of all, we prove that the remainder converge in probability with a $1/\sqrt{n}$ convergence rate.
 On account of Theorem \ref{Orcalebounds}, we get
\begin{equation}
\begin{split}
||rem||_{\infty}\leq \sum_{k=1}^{2} ||(\widehat{\Theta}^{[k]}-\Theta_0^{[k]})
 (\widehat{\Sigma}^{[k]}-\Sigma_0^{[k]})\Theta_0^{[k]}||_{\infty}
 +\sum_{k=1}^{2} || (\widehat{\Theta}^{[k]}-\Theta_0^{[k]})
 (\widehat{\Sigma}^{[k]}\widehat{\Theta}^{[k]}-\bbI_p) ||_{\infty}
\end{split}
\end{equation}
 Define
\bqa
l(\Theta)=\sum_{k=1}^{2} \{tr(\widehat\Sigma^{[k]}\Theta^{[k]})-\log \det(\Theta^{[k]})\}+
\lambda\sum_{k=1}^{2}||(\Theta^{[k]})^{-}||_{1}
+\rho||(\Theta^{[1]}
-\Theta^{[2]})^{-}||_1.
\eqa
 By the Karush-Kuhn-Tucker (KKT) conditions, we yield
\bqa\label{KKT}
\widehat\Sigma^{[1]}-(\widehat{\Theta}^{[1]})^{-1}+(\lambda+\rho)\widehat{Z}^{[1]}=0,
\eqa
 and
\bqa
\widehat\Sigma^{[2]}-(\widehat{\Theta}^{[2]})^{-1}+(\lambda-\rho)\widehat{Z}^{[2]}=0,
\eqa
 where $\widehat{Z}^{[k]}_{ij}=sign(\widehat{\Theta}^{[k]}_{ij})$ if $\widehat{\Theta}^{[k]}_{ij}\neq 0$, and satisfying $||  \widehat{Z}^{[k]} ||_{\infty}\leq 1$.
 Multiplying by $\widehat{\Theta}^{[1]}$ on the equation (\ref{KKT}), we get
\bqa
\bbI_p-\widehat\Sigma^{[1]}\widehat{\Theta}^{[1]}
=(\lambda+\rho)\widehat{Z}^{[1]}\widehat{\Theta}^{[1]}.
\eqa
 Similarly, we have
\bqa
\bbI_p-\widehat\Sigma^{[2]}\widehat{\Theta}^{[2]}
=(\lambda-\rho)\widehat{Z}^{[2]}\widehat{\Theta}^{[2]}.
\eqa
 Thus,
\begin{equation}
\begin{split}
||rem||_{\infty}\leq& \sum_{k=1}^{2} |||  (\widehat{\Theta}^{[k]}-\Theta_0^{[k]})  |||_1
|| (\widehat{\Sigma}^{[k]}-\Sigma_0^{[k]})  ||_{\infty}
|||  \Theta_0^{[k]} |||_1 \\
&+(\lambda+\rho)\sum_{k=1}^{2} |||  (\widehat{\Theta}^{[k]}-\Theta_0^{[k]}) |||_1
||  \widehat{Z}^{[k]} ||_{\infty}
|||   \widehat{\Theta}^{[k]}  |||_1.
\end{split}
\end{equation}
 To draw the conclusion, we have
\bqa\label{colsumnormb}
|||  (\widehat{\Theta}^{[k]}-\Theta_0^{[k]})  |||_1 \leq b(p+s)\lambda,
\eqa
 where $b$ is a constant and is related to $L$.
 According to the Schwarz inequality and Weyl inequality, we get
\bqa
|||  \Theta_0^{[k]} |||_1\leq \sqrt{d+1}\Lambda_{\max}(\Theta_0^{[k]}).
\eqa
 The bound of $|||   \widehat{\Theta}^{[k]}  |||_1$ is derived by
\bqa
|||   \widehat{\Theta}^{[k]}  |||_1\leq|||   \widehat{\Theta}^{[k]} - \Theta_0^{[k]}|||_1+|||\Theta_0^{[k]}|||_1.
\eqa
 According to the rate of $\lambda$, we conclude that
\bqa\label{colsumnorme}
|||   \widehat{\Theta}^{[k]}  |||_1\leq \sqrt{d+1}\Lambda_{\max}(\Theta_0^{[k]}).
\eqa
 Besides, the Sub-Gaussian random vector with covariance $\Sigma_0^{[k]}$ implies that $||\widehat\Sigma^{[k]}-\Sigma_0^{[k]}||_{\infty}=O_p(\sqrt{\log(p)/n})$, where $O_p$ denotes bounded in probability.
 We get
\bqa
||rem||_{\infty}&\leq&\frac{4\lambda(8s_1+8s_2+p)}{c}\sqrt{\frac{\log{p}}{n}}
\sqrt{d+1}\max\{\Lambda_{\max}(\Theta_0^{[1]}),\Lambda_{\max}(\Theta_0^{[2]})\}\ \non
&&+(\lambda+\rho)\frac{4\lambda(8s_1+8s_2+p)}{c}\sqrt{d+1}
\max\{\Lambda_{\max}(\Theta_0^{[1]}),\Lambda_{\max}(\Theta_0^{[2]})\}.
\eqa
 For $\lambda\asymp \rho$, $||rem||_{\infty}$ is bounded by $\tilde{b}(p+s)\sqrt{d+1}\lambda^2$ in probability, where $\tilde b$ is a constant related to $L$.
 Based on the condition $(p+s)\sqrt{d}=o(\sqrt{n}/\log{p})$, $||rem||_{\infty}=o_{p}(1/\sqrt{n})$.
 According to the bounded fourth moments of $(\widehat{\Theta}^{[k]})_{ii}(\widehat{\Theta}^{[k]})_{jj}
 +(\widehat{\Theta}^{[k]})_{ij}^2$ and Lindeberg central limit theorem, we complete the proof of the Theorem \ref{CLT}.

\end{proof}

\subsection{Proof of Theorem \ref{weiverclt}}

\begin{proof}
 The conclusions of Theorem \ref{weiverclt} can be obtained from the arguments (\ref{colsumnormb})-(\ref{colsumnorme}). For weighted version, $||rem||_{\infty}$ can be bounded by $\tilde{b}s\sqrt{d+1}\lambda^2$, which completes the proof.
\end{proof}

\end{appendices}

\end{document}